\documentclass[preprint,12pt]{elsarticle}

\usepackage{amsmath,amssymb,graphicx, float, algorithm2e, booktabs, varwidth, color, hyperref, adjustbox, multirow, multicol, centernot, mathrsfs, pifont, tikz, enumerate}
\usetikzlibrary{hobby}
\usetikzlibrary{patterns}
\usetikzlibrary{calc}

\newcommand\ellipsebyfoci[4]{
  \path[#1] let \p1=(#2), \p2=(#3), \p3=($(\p1)!.5!(\p2)$)
  in \pgfextra{
    \pgfmathsetmacro{\angle}{atan2(\y2-\y1,\x2-\x1)}
    \pgfmathsetmacro{\focal}{veclen(\x2-\x1,\y2-\y1)/2/1cm}
    \pgfmathsetmacro{\lentotcm}{\focal*2*#4}
    \pgfmathsetmacro{\axeone}{(\lentotcm - 2 * \focal)/2+\focal}
    \pgfmathsetmacro{\axetwo}{sqrt((\lentotcm/2)*(\lentotcm/2)-\focal*\focal}
  }
  (\p3) ellipse[x radius=\axeone cm,y radius=\axetwo cm, rotate=\angle];
}
\SetKw{Set}{set}{}
\SetKwBlock{Repeat}{repeat}{}
\newcommand{\ee}[1]{\times10^{#1}} 
\newcommand{\norm}[1]{\Vert #1 \Vert} 
\newcommand{\eps}{\varepsilon} 

\newtheorem{defn}{Definition}
\newcommand{\aref}[1]{Algorithm \ref{#1}}
\newcommand{\eref}[1]{(\ref{#1})}
\newcommand{\Fref}[1]{Figure~\ref{#1})}

\journal{Applied Numerical Mathematics}

\begin{document}
\begin{frontmatter}
\title{Plug-and-play superiorization}

\author[uwb]{Jon Henshaw}
\author[hit]{Aviv Gibali}
\ead{avivgi@hit.ac.il}
\author[uwb]{Thomas Humphries\corref{cor1}}
\ead{thumphri@uw.edu}

\affiliation[uwb]{organization={Department of Engineering and Mathematics, School of STEM, University of Washington Bothell},
city={Bothell},
postcode={98011},
state={WA},
country={USA}}

\affiliation[hit]{organization={Applied Mathematics Department, Holon Institute of Technology},
city={Holon},
postcode={5810201},
country={Israel}}

\cortext[cor1]{Corresponding author}

\begin{abstract} The superiorization methodology (SM) is an optimization heuristic in which an iterative algorithm, which aims to solve a particular problem, is ``superiorized'' to promote solutions that are improved with respect to some secondary criterion. This superiorization is achieved by perturbing iterates of the algorithm in nonascending directions of a prescribed function that penalizes undesirable characteristics in the solution; the solution produced by the superiorized algorithm should therefore be improved with respect to the value of this function. In this paper, we broaden the SM to allow for the perturbations to be introduced by an arbitrary procedure instead, using a plug-and-play approach. This allows for operations such as image denoisers or deep neural networks, which have applications to a broad class of problems, to be incorporated within the superiorization methodology. As proof of concept, we perform numerical simulations involving low-dose and sparse-view computed tomography image reconstruction, comparing the plug-and-play approach to two conventionally superiorized algorithms, as well as a post-processing approach. The plug-and-play approach provides comparable or better image quality in most cases, while also providing advantages in terms of computing time, and data fidelity of the solutions. 

\end{abstract}


\end{frontmatter}



\section{Introduction}

Many real-world problems in science and technology can be modeled as optimization problems featuring several types of constraints. For example, in image reconstruction, one often wishes to find an image that fits some measured data, and also satisfies some secondary criteria based on prior knowledge (e.g. sparsity in some representation, or smoothness). As it is often not possible -- or desirable, in the presence of noise -- to exactly fit the measured data, there exist many images that are feasible in the sense of attaining the same level of data fidelity (as measured by some proximity function), but vary in quality with respect to the secondary criteria. Achieving both feasibility and {\em optimality} with respect to the secondary criteria is often very demanding computationally. If one does not need the best solution, however, but only a very good one, it is beneficial to resort to suboptimal but still feasible solutions. This is exactly the idea behind the {\em superiorization methodology} (SM) on which this work is built: take a feasibility seeking algorithm, and develop a framework to steer the outcomes of this algorithm towards better solutions, while keeping the computational effort low.

The idea of the SM is to start from an initial solution of the problem and then to introduce certain perturbations to the solution such that the perturbed result is superior to the original solution in some well defined sense, but not necessarily optimal. By trading optimality for superiority, one has more flexibility in designing new methods, and can achieve a considerable algorithmic speedup. Superiorization has successfully found multiple applications, such as computed tomography \cite{cdhst14,gcohj22}, inverse treatment planning in radiation therapy \cite{dcsgx15}, bioluminiscence tomography \cite{jcj13} and linear optimization \cite{c17}, in some cases outperforming other state-of-the-art algorithms. The state of current research on superiorization is tracked on a continually updated bibliographical website \cite{CensorSupPage}; special issues of {\em Inverse Problems} and {\em Journal of Applied and Numerical Optimization} focused on the SM have also been published~\cite{CHJ17,GHS20}.

A key feature of the SM is its flexibility; that is, given one's preferred feasibility seeking algorithm (the {\em basic algorithm}), superiorization provides an automatic procedure to generate improved solutions, provided that the basic algorithm is {\em perturbation resilient}. To date, all applications of superiorization define a penalty function $\phi$ and perturb iterates in a nonascending direction of $\phi$, thus promoting superiority of the solution with respect to this function. In many applications, however, a user may have some other procedure which improves the quality of a solution; for example, a denoiser when problems involve image restoration or reconstruction. Additionally, there has been a recent explosion of interest in methods based on deep learning in the field of image reconstruction, replacing model-driven penalty functions with data-driven prior information encoded within a deep neural network. In both of these cases, the procedure may be viewed as a ``black box'' where the underlying function being minimized is intractable, making it difficult to apply within a framework requiring nonascending directions.

Motivated by these observations, we propose a ``plug-and-play'' superiorization approach which allows one to introduce perturbations generated by an arbitrary procedure within the SM. The idea is that the procedure in question improves the quality of the current iterate with respect to some criterion, without explicitly reducing the value of some target function. This approach is inspired by the success of the plug-and-play alternating direction method of multipliers (ADMM)\cite{vbw13} on a variety of problems, including image restoration~\cite{cwe17} and low-dose CT reconstruction~\cite{hywzb18}. 

The paper is organized as follows. In Section~\ref{S:method}, we introduce the basic concepts behind the conventional SM, as well as our proposed plug-and-play methodology. Additionally, we describe the mathematical model used for the CT image reconstruction problem considered in this paper, as well as two challenging scenarios within the field. Section~\ref{S:exp} describes our numerical experiments in which we study the effectiveness of the plug-and-play approach, using both a state-of-the-art denoising algorithm, as well as a neural network to generate perturbations; results are presented in Section~\ref{S:results}. We include some discussion and concluding remarks in Sections~\ref{S:disc} and \ref{S:concl}.

\section{Methodology}\label{S:method}
\subsection{Superiorization methodology}\label{S:method_sup}
The superiorization methodology (SM) is based on the following concepts.  We define a problem set $\mathbb{T}$ to be a general description of the problem to be solved, and $T \in \mathbb{T}$ to be a specific problem from this set. For example, $\mathbb{T}$ could describe a general convex feasibility problem

\begin{eqnarray}
\textrm{Find } x^* \in \bigcap_{i = 1}^I C_i \label{E:cfprob}
\end{eqnarray}
where the $C_i$ are convex sets in $\mathbb{R}^J$; a problem $T$ would then specify the particular $C_i$ to whose intersection we wish to converge. We then specify an algorithm for problem set $\mathbb{T}$ which defines, for each $T \in \mathbb{T}$, an operator $\boldsymbol{P}_{T}:\mathbb{R}^J \to \mathbb{R}^J$. From any $x^0 \in \mathbb{R}^J$, this operator, called the {\em basic algorithm}, produces an infinite sequence via the iteration
\begin{eqnarray}
x^{k+1} = \boldsymbol{P}_{T}(x^k) = \boldsymbol{P}_T^{k+1}(x^0). \label{E:basic}
\end{eqnarray}
We then have the following~\cite{cdh10}:

\begin{defn}
A basic algorithm is {\bf bounded perturbation resilient} if, whenever \eref{E:basic} converges to a solution of the problem $T$ for all $x^0\in\mathbb{R}^{J}$, then so does the iteration
\begin{eqnarray}
x^{k+1} = \boldsymbol{P}_{T}(x^k+\beta_{k}\nu
^{k}), \label{E:sup}
\end{eqnarray} where $\beta_{k}$ are real nonnegative numbers such that $\sum_{i=0}^{\infty}\beta_{k}<\infty$, and the sequence $\{\nu^{k}\}_{k=0}^{\infty}$ is bounded. \label{D:bpr}
\end{defn}

In many contexts one may not wish to exactly solve the problem $T$; for example, if the constraints to be satisfied involve noisy data. In this case the related notion of \textbf{strong
perturbation resilience}~\cite{cdhst14,hgdc12} is useful. We define a {\em proximity function} $Pr_T:\mathbb{R}^J \to \mathbb{R}_+$ which indicates the extent to which a given $x \in \mathbb{R}^J$ satisfies the constraints of the problem; for example, the sum of Euclidean distances of $x$ from each $C_i$ in \eref{E:cfprob}. We say that $x$ is \textit{$\varepsilon$-compatible} with $T$ if $Pr_T(x) \leq \varepsilon$. Finally, given a particular $\varepsilon > 0$, we define the {\em output} of iteration~\eref{E:basic} to be the first $x^k$ that is $\varepsilon$-compatible, if such an iterate exists (otherwise, the output is undefined for that value of $\varepsilon$). We then have the following~\cite{hgdc12}:

\begin{defn}
A basic algorithm is {\bf strongly perturbation resilient} if, for all $T \in \mathbb{T}$,
\begin{enumerate}
    \item There exists an $\varepsilon > 0$ such that the output of~\eref{E:basic} is defined for all $x^0 \in \mathbb{R}^J$, and
    \item For all $\varepsilon > 0$ such that the output of~\eref{E:basic} is defined for all $x \in \mathbb{R}^J$, the output of \eref{E:sup} is also defined for all $x^0$ and all $\varepsilon' > \varepsilon$.
\end{enumerate}\label{D:spr}
\end{defn}

\begin{figure}[H]
    \centering
    \begin{adjustbox}{width=\linewidth}
    \begin{tikzpicture}
        \ellipsebyfoci{fill, pattern=north east lines, pattern color=green!50!black, use Hobby shortcut, closed=true}{1.5,3.75}{3.5,4.5}{1.5}
        \ellipsebyfoci{draw,color=green!50!black}{1.5,3.75}{3.5,4.5}{1.5}
        \node[color=green!50!black, fill=green!10!white] (c) at (2.55,4.305) {$Pr_T(x)\leq\eps$};

        \path[draw,color=black,-stealth,style=ultra thick,use Hobby shortcut,closed=false] (-.25,3) .. (1.4, 3.54);

        \filldraw[black] (-.25,3) circle (2pt);

        \path[draw,color=black,-stealth,style=ultra thick,use Hobby shortcut,closed=false] (-1.5, 5) .. (0.75,3.75);
        
        \node[] (c) at (0,3.75) {$x^{k+1}$};
        
        \path[draw,color=black,-stealth,style=ultra thick,use Hobby shortcut,closed=false] (0.75,3.75) .. (-.25,3);
        
        \filldraw[black] (0.75,3.75) circle (2pt);

        \path[draw,color=red,-stealth,style=ultra thick,use Hobby shortcut,closed=false] (-0.75, 6.5) .. ((1.3, 4.875);
        
        \path[draw,color=black,-stealth,style=ultra thick,use Hobby shortcut,closed=false] (-0.75, 6.5) .. (-1.5, 5);
        
        \filldraw [black] (-1.5, 5) circle (2pt);

        \filldraw [black] (-0.75, 6.5) circle (2pt);

        \node[] (c) at (-1.07,6.85) {$x^k$};
        \node[] (c) at (2,4.9) {$x^{k+1}$};
        \node[] (c) at (-2.5,5.2) {${x}^{k}+\beta_k\nu^k$};
        \node[] (c) at (-1.5,3.10) {${x}^{k+1}+\beta_{k+1}\nu^{k+1}$};
        \node[] (c) at (2,3.6) {${x}^{k+2}$};

        \filldraw [color=black,fill=red!90!black] (1.406845567387103, 4.84303310507534) circle (2.5pt);

        \filldraw [color=black,fill=green!90!black] (1.4, 3.54) circle (2.5pt);
        
        \filldraw [color=black,fill=blue!90!black] (-2.5, 2) circle (2.5pt);
        \node[] (c) at (-3.5, 1.7) {$\displaystyle \arg \min_x ~\phi$};
        \draw[color=blue!90!black] (-1.9,1.2) arc
            [   color=blue,
                start angle=-45,
                end angle=135,
                x radius=1cm,
                y radius =1cm
            ] ;
        \draw[color=blue!90!black] (-1,1) arc
            [   color=blue,
                start angle=-45,
                end angle=135,
                x radius=2cm,
                y radius =2cm
            ] ;
        \draw[color=blue!90!black] (0,0) arc
            [   color=blue,
                start angle=-45,
                end angle=135,
                x radius=3.5cm,
                y radius =3.5cm
            ] ;
    \end{tikzpicture}
    \end{adjustbox}
    \caption{Illustration of the superiorization process. The iterate $x^{k+1}$ attained by the basic algorithm (red circle) defines a region of $\varepsilon$-compatibility (green shaded ellipse). In the superiorized algorithm (black arrows), the solution is perturbed in non-ascending directions of the function $\phi$, whose level curves are shown in blue. Eventually, the superiorized solution $x^{k+2}$ (green circle) is attained. Note that this solution is not necessarily the optimal $\varepsilon$-compatible solution.}\label{F:sup}
    \label{blob}
\end{figure}
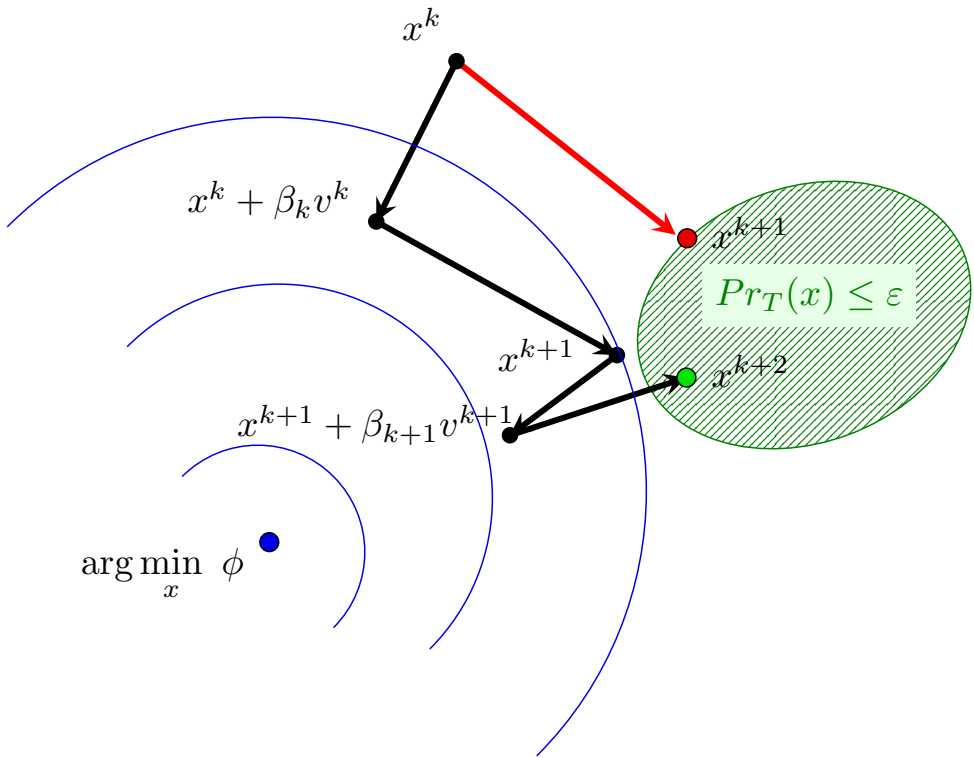

Loosely speaking, if an algorithm is strongly perturbation resilient, then whenever \eref{E:basic} converges to an $\varepsilon$-compatible solution, \eref{E:sup} will eventually do so as well. The essential idea of the superiorization methodology is to make use of the perturbations to transform a bounded (or strongly) perturbation-resilient algorithm into an algorithm whose outputs are equally good from the point of view of constraints compatibility, but are superior (but not necessarily optimal) according to some other criteria, typically measured by a penalty function $\phi : \mathbb{R}^J \to \mathbb{R}$. This is done by producing from the basic algorithm another algorithm, called its \textit{superiorized} version, that makes sure not only that the $\beta
_{k}\nu^{k}$ are bounded perturbations, but also that $\phi\left(  x^{k}+\beta_{k}\nu^{k}\right)  \leq\phi\left(  x^{k}\right) $, for all $k\geq 0$. This is achieved by choosing the $\nu^k$ to be nonascending directions of the function $\phi$ at $x^k$. \Fref{F:sup} provides an illustration of the process.


We present a typical framework for a superiorized algorithm in \aref{A:sup_conv}. For concreteness, we make the following decisions and assumptions:
\begin{enumerate}[(i)]
    \item To promote descent of the penalty function, we introduce multiple perturbations in a nonascending direction of $\phi$ between each iteration of the basic algorithm. The parameter $N \in \mathbb{N}$ specifies the number of such perturbations. This introduces a second loop, such that each perturbation is of the form $x^{k,n} + \beta_{k,n} \nu^{k,n}$, for $n \in [0, N)$, with $x^{k+1} = \mathbf{P}_T(x^{k,N})$.
    \item To ensure summability of the $\beta_{k,n}$, we introduce parameters $\alpha > 0$, $\gamma \in (0,1)$, such that $\beta_{0,0} = \alpha$, and subsequent values decrease geometrically with ratio $\gamma$.
    \item We assume that $\phi$ is differentiable everywhere, so that $\nu^{k,n}$ can simply be chosen as the negative gradient of $\phi$ at $x^{k,n}$, normalized to unit length.
\end{enumerate}

Both (i) and (ii) are standard choices in the SM literature. In the case where $\phi$ is not everywhere differentiable, subgradients or more generic choices of nonascending direction can be used instead, see e.g.~\cite{CZ15,cghh21,ACGT23}.

\begin{algorithm}[H]
\caption{Conventionally superiorized version of a basic algorithm}
    \label{A:sup_conv}
Given $x^{0}$, $\varepsilon > 0$, $N \in \mathbb{N}$, $\gamma \in (0, 1)$, $\alpha > 0$: \\
\medskip
$\ell \gets -1$ \\
\For{$k = 0, 1, 2, \dots$}{
    $x^{k,0} \gets x^{k}$ \\
    \For{$n = 0, 1, 2,  \dots N-1$}{
        $\nu^{k,n} \gets -\nabla \phi (x^{k,n}) / \left( \norm{ \nabla \phi(x^{k,n}}\right)$ \\
        \While{true}{
            $\ell \gets \ell+1$ \\
            $\beta_{k,n} \gets \alpha \gamma^\ell$ \\
            $z \gets x^{k,n} + \beta_{k,n} \nu^{k,n}$ \\
            \If {$\phi(z) < \phi(x^{k})$}{
                 $x^{k,n+1} \gets z$ \\
                {\bf break} \\
            }
        }
    }
    $x^{k+1} \gets  \mathbf{P}_T (x^{k,N})$ \\
    \If {$Pr_T(x^{k+1}) < \varepsilon$}{
        \Return $x^{k+1}$ \\
    }
}
\end{algorithm}

The performance of a superiorized algorithm is often sensitive to the parameters $N$, $\gamma$, $\varepsilon$. In general, choosing $N$ to be larger puts more effort towards reducing the value of $\phi(x^k)$ between every iterationn, as more steps in nonascending direction of $\phi$ are taken; similarly, choosing $\gamma$ close to 1 causes the size of the perturbations to decay slowly, resulting in larger steps taken in nonascending directions of $\phi$ as the algorithm proceeds.  Choosing either value to be too small can result in the superiorized algorithm producing solutions that are not substantially better with respect to $\phi$ than those produced by the basic algorithm. On the other hand, choosing them to be too large results in slow convergence, and can also result in poor solutions. Similarly, the choice of stopping criterion $\varepsilon$ is important; if the value is too large, the algorithm may not yet have converged to an accurate solution, while if it is too small (resulting in a need for more iterations), the superiorized algorithm may produce solutions similar to those of the basic algorithm, because the perturbation size becomes very small in later iterations. In practice, good values of these parameters vary depending on the problem being solved, and are often chosen after some experimentation; see~\cite{humphries2020comparison} for some discussion.

Recently, a new rule for selecting the step size for the SM was proposed in~\cite{naae22}. The new algorithm eliminates the loop over $N$ in \aref{A:sup_conv}, and has the form:

\begin{algorithm}[H]
\caption{Superiorized algorithm with adaptive step size}
    \label{A:sup_adap}
Given $x^{0}$, $\varepsilon > 0$, parameters $\alpha_0, \epsilon$: \\
\medskip
\For{$k = 0, 1, 2, \dots$}{
    $\nu^{k} \gets -\nabla \phi (x^{k}) / \left( \norm{ \nabla \phi(x^{k}}\right)$ \\
   \uIf {$\phi(x^k) < \alpha_k$}{
        $\beta_k \gets 0$ \\
    }\Else{
        $\beta_k \gets \frac{\phi(x^k)-\alpha_k}{\norm{ \nabla \phi(x^{k}}}$ \\
    }
    $z \gets x^{k} + \beta_{k} \nu^{k}$ \\  
    $\zeta_k \gets \left( Pr_T(z) - Pr_T(x^k)\right) / Pr_T(x^k)$\\
    $\alpha_{k+1} \gets \alpha_k + \max\{ \epsilon,-\zeta_k \alpha_k\}$\\
    $x^{k+1} \gets  \mathbf{P}_T (z)$ \\
    \If {$Pr_T(x^{k+1}) < \varepsilon$}{
        \Return $x^{k+1}$ \\
    }
}
\end{algorithm}
In this approach, the superiorization step functions as a subgradient projector relative to $\phi$ by level $\alpha_k$~(\cite{C12}, Definition 4.2.4). The quantity $\zeta_k$ is called the {\em desirability number} and controls the increase of the sequence of levels $\{\alpha_k\}$. The idea is that if $\zeta_k > 0$, reducing the value of $\phi(x^k)$ below $\alpha_k$ has caused the proximity function $Pr_T$ to increase, and vice-versa if $\zeta_k < 0$. When the data are noisy, we wish to prioritize reducing $\phi$ over reducing $Pr_T$, and so positive values of $\zeta_k$ are desirable; in this case, $\{\alpha_k\}$ increases slowly (by the parameter $\epsilon$), improving effectiveness of the method. The reverse is true when the data are noiseless, and so the updating rule $\alpha_{k+1} \gets \alpha_k + \max\{ \epsilon,\zeta_k \alpha_k\}$ is used in this case~\cite{naae22}.


\subsection{Plug-and-Play superiorization}

Suppose now that $\Psi: \mathbb{R}^J \to \mathbb{R}^J$ is a procedure designed to improve the quality of a given $x$ in some way; for example, a denoiser, or a deep neural network trained to remove other types of artifact from an image. If we consider applying this procedure in between iterations of the basic algorithm, we obtain the iteration
\begin{eqnarray}
x^{k+1}  &= \boldsymbol{P}_T( \Psi (x^k)) = \boldsymbol{P}_T( x^k + v^k )\label{E:pnp1}
\end{eqnarray}
where $v^k = \Psi(x^k) - x^k$; i.e., the action of $\Psi$ implicitly defines a perturbation to $x^k$. Comparing with \eref{E:sup}, we can see a straightforward analogy to superiorization, with the critical exception that \eref{E:pnp1} offers no guarantee that the perturbations are bounded, let alone that their norms are summable. We can, however, define a modified iteration with exactly the same general form as \eref{E:sup} by defining
\begin{eqnarray}
    v^k &= \Psi(x^k) - x^k, \\
    \nu^k &= v^k / \Vert v^k\Vert, \\
    \beta_k &= \min \left\{ \mu_k, \Vert v^k \Vert \right \},
\end{eqnarray}
with $\sum_{k=0}^\infty \mu_k < \infty$. As before, we choose $\mu_k = \alpha \gamma ^k$ for $\gamma \in (0,1)$ to ensure that this holds. The result is that if the perturbation induced by $\Psi$ is sufficiently small, the basic algorithm is simply applied to $\Psi(x^k)$; otherwise, the perturbation $v^k$ is damped by a factor $\mu_k / \Vert v^k \Vert$ that ensures summability.

The motivation for the plug-and-play superiorization (PnP-sup) approach is twofold:

\begin{enumerate}
    \item While this represents a straightforward extension of the existing SM -- in that no new theory is required -- it allows for considerably more flexibility in the design and implementation of superiorized algorithms. For example, in the case where $\Psi$ is a neural network, the same network architecture could be trained on different datasets to achieve different goals, such as denoising, artifact removal, or even acceleration of the basic iteration (by training it to map early iterates to ones resembling later iterates). This avoids the need to construct different model-based priors to achieve these goals. As discussed in Section~\ref{S:exp:nn}, the use of neural networks has attracted particular interest in the field of tomographic imaging in recent years. Given that tomography is one of the main use cases of the SM, dating back to some of the earliest papers~(e.g.~\cite{cdh10,hgdc12}), the ability to incorporate neural networks within the SM is especially enticing.
    \item Procedures such as denoisers, neural networks, etc. are often applied to a solution generated by the basic algorithm as a post-processing step. While this may produce a solution with desirable characteristics, there is no guarantee that the output of this post-processing step is $\varepsilon$-compatible, even if the output of the basic algorithm was. The PnP-sup methodology ensures $\varepsilon$-compatibility while also making use of the procedure to improve solution quality within the iteration.
\end{enumerate}

The overall algorithm is presented below as \aref{A:sup_pnp}. Convergence of the algorithm to a $\varepsilon$-compatible solution follows directly from the definition of strong perturbation resilience (Definition~\ref{D:spr}) and the fact that $\norm{\nu_k}$ = 1, $\sum \beta_k < \sum \alpha \gamma^k < \infty$, meeting the requirements on the sequence of iterates given in Definition~\ref{D:bpr}. This is somewhat more straightforward than the proof of convergence of \aref{A:sup_conv} presented in~\cite{hgdc12},  which requires aggregating the $N$ perturbations in the inner loop into a single bounded perturbation in order to apply the definition.

\begin{algorithm}[H]
\caption{Plug-and-play superiorization (PnP-sup)}
    \label{A:sup_pnp}
Given $x^{0}$, $\varepsilon > 0$, $\gamma \in (0, 1)$, $\alpha > 0$: \\
\medskip
$\ell \gets -1$ \\
\For{$k = 0, 1, 2, \dots$}{
    $z \gets \Psi(x^{k})$ \\
    $v^k \gets z - x^{k}$ \\
    $\ell \gets \ell + 1$\\
    $\beta_k \gets \min \{\alpha \gamma^\ell, \norm{v^k} \}$ \\
    $x_+^k \gets x^k + \frac{\beta_k}{\norm{v^k}} v^k$\\
    $x^{k+1} \gets  \mathbf{P}_T (x_+^k)$ \\
    \If {$Pr_T(x^{k+1}) < \varepsilon$}{
        \Return $x^{k+1}$ \\
    }
}
\end{algorithm}

\subsection{Computed tomography}

\begin{figure}
    \centering
    \includegraphics[width=0.53\linewidth]{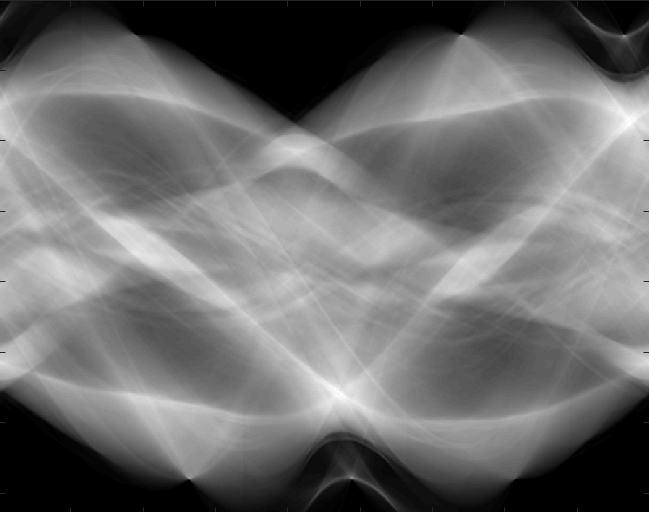}
    \includegraphics[width=0.42\linewidth]{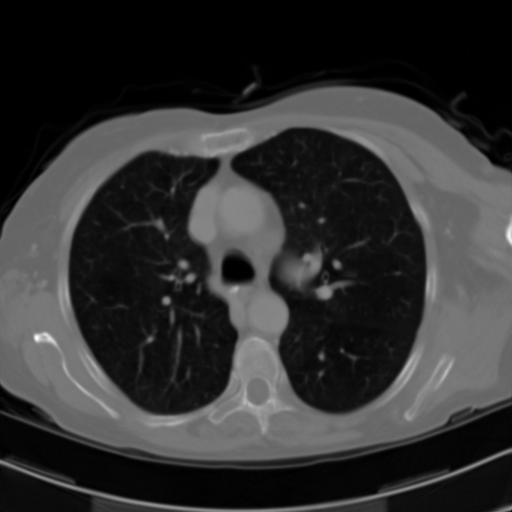}
    \caption{A simulated fan-beam sinogram (left) and the corresponding reconstructed CT image slice (right). Each column of the sinogram represents fan-beam X-ray measurements taken from a particular source position. Image courtesy of The Cancer Imaging Archive~\protect\cite{TCIA,QIN_LUNG_CT}, with sinogram produced using the ASTRA toolbox~\protect\cite{PBS11,VPCJB16}. }
    \label{F:sino}
\end{figure}

Computed tomography (CT) is a medical imaging technique that allows for the reconstruction of cross-sectional or volumetric images of the body via a series of X-ray measurements through the patient. Given a monoenergetic X-ray beam of intensity $I_0$ counts per unit time, and a spatially varying, compactly supported attenuation function $\mu(x)$ in two or three dimensions, the Beer-Lambert law states that the idealized X-ray measurement along a line $j$ is given by
\begin{eqnarray}
\hat{I}_j &&= I_0 \exp \left( - \int_j \mu(x) \: dx \right), \textrm{ or } \label{E:beers}\\
\ln( I_0 / \hat{I}_j ) &&= \int_j \mu(x) \: dx. \nonumber
\end{eqnarray}
Thus, the standard mathematical model is that the log-transformed X-ray measurements represent line integrals through the distribution of attenuating tissue within the body, i.e., samples of the X-ray transform (which is equivalent to the Radon transform for functions of two variables). Discretizing in both the image and measurement space leads to a large linear system of equations:
\begin{eqnarray}
    Ax = b \label{E:linear}
\end{eqnarray}
where $x \in \mathbb{R}^J$ represents the pixelized (or voxelized) image to be reconstructed, $b \in \mathbb{R}^I$ is the measured data (referred to as the sinogram), and $A$ is the $I \times J$ system matrix whose $(i,j)$th entry represents the contribution of the $j$th image pixel to the $i$th measurement. Note that \eref{E:linear} can be viewed as a special case of problem \eref{E:cfprob}, where $C_i = \{ x \mid a_ix = b_i\}$, with $a_i$ denoting the $i$th row of matrix $A$. A 2D CT image slice $x$ and its corresponding simulated sinogram $b$ are shown in \Fref{F:sino}

In recent years there has been considerable interest in reducing the patient dose associated with CT scans. Two ways of accomplishing this goal are:
\begin{enumerate}
    \item Reducing the beam intensity, given as $I_0$ in \eref{E:beers}. The main drawback of this approach (which we refer to as low-dose CT) is that the measured intensity of a beam traveling along line $j$, $I_j$, is a random variable, due to the stochastic nature of the interaction of X-rays with matter. This is typically modeled as having a Poisson distribution with mean $\hat{I}_j$ and standard deviation (and therefore, signal-to-noise ratio) of $\sqrt{\hat{I}_j}$. Thus, a reduction in $I_0$ causes a corresponding reduction in $\hat{I}_j$, producing noisier data and a noisier reconstructed image.
    \item Reducing the number of angular samples, typically known as sparse-view CT. In this approach, the X-ray source is rapidly switched on and off to acquire data from fewer source positions than are typically used. The resulting undersampling produces streaking artifacts in images reconstructed using conventional algorithms.
\end{enumerate}

\begin{figure}
    \centering
    \includegraphics[width=0.42\linewidth]{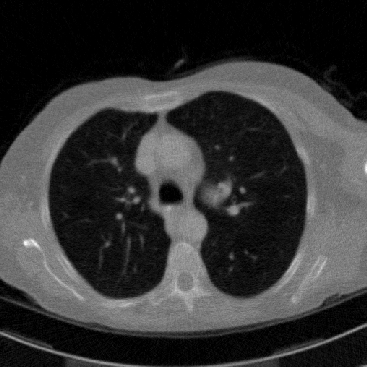}
    \includegraphics[width=0.42\linewidth]{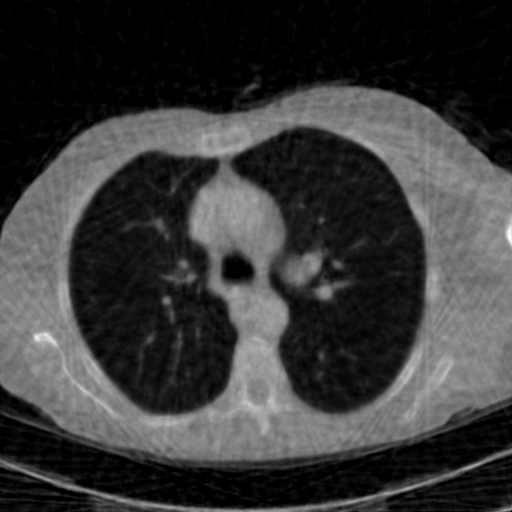}
    \caption{CT slices reconstructed from low-dose data (left) and sparse-view data (right). Compare with \Fref{F:sino}, right hand side. }\label{F:lowdose}
\end{figure}

Examples of images produced from low-dose and sparse-view measured data are shown in \Fref{F:lowdose}. Both of these scenarios have been studied extensively in recent years, using mathematical techniques such as total variation minimization or other compressive-sensing based approaches, or, more recently, with deep neural networks. We consider both scenarios in this paper to demonstrate the effectiveness of the PnP-sup approach.

\subsection{Iterative reconstruction}

A number of algorithms have been used to solve \eref{E:linear} iteratively in the context of CT imaging. In this work we use a variant of the Simultaneous Iterative Reconstruction Technique (SIRT), which has the form

\begin{eqnarray}
    x^{k+1} &= x^k - \omega_k D A^T M(Ax^k-b)\label{E:SART},
\end{eqnarray}
where $\omega_k \in (0,2)$ is a relaxation parameter, and $D,M$ are $J \times J$ and $I \times I$ diagonal matrices, respectively, whose entries are given by:

\begin{eqnarray}
    D_{jj} &= 1 \bigg/\sum_{k=1}^{I} |a_{kj}|,~~M_{ii} &=  1 \bigg/\sum_{k=1}^{J} |a_{ik}|; \label{E:DandM}
\end{eqnarray}
i.e. the reciprocals of the column and row sums of $A$. SIRT tends to converge slowly, but can be accelerated by partitioning the sinogram into ordered subsets, consisting of equally spaced views of the object being imaged. For example, if ten subsets are used, the first consists of the 1st, 11th, 21st, etc. views, the second of the 2nd, 12th, 22nd, and so on. Let $s(w)$ denote the $w$th subset of projection data, $w = 1 \dots W$, and let $A_{s(w)}$ be the matrix obtained by extracting only the rows of the system matrix $A$ corresponding to $s(w)$, and similarly for $b_{s(w)}$, $D_{s(w)}$ and $M_{s(w)}$. Then, the basic algorithm we use in this paper is given by
\begin{eqnarray}
\mathbf{P}_T(x) = Proj_{\mathbb{R}_+} \mathbf{B}_{W} \dots \mathbf{B}_2 \mathbf{B}_1 (x), \label{E:OS-SIRT1}
\end{eqnarray}
where
\begin{eqnarray}
\mathbf{B}_w (x) = x  - \omega_k D_{s(w)} (A_{s(w)})^T M_{s(w)}\left[ A_{s(w)} x - b_{s(w)} \right],\label{E:OS-SIRT2}
\end{eqnarray}
and $Proj_{\mathbb{R}_+}$ is the projection onto the non-negative orthant, which sets any negative entries of $x$ to zero. This approach has been referred to as OS-SIRT~\cite{BKK12}; the well-known simultaneous algebraic reconstruction technique (SART)~\cite{andersen1984simultaneous} corresponds to the case where each subset consists of a single view. We refer to it as block-iterative SART (BI-SART), consistent with~\cite{CE02}.

In \cite{hgdc12} it is proved that an algorithm $\boldsymbol{P}_T$ is strongly perturbation resilient (Definition~\ref{D:spr}) provided that it is boundedly convergent and non-expansive, and that the proximity function used in the definition, $Pr_T(x)$, is uniformly continuous. It is known that the SIRT algorithm~\eref{E:SART} converges to a solution to the weighted least-squares problem
\begin{eqnarray}
    \min_x \frac{1}{2} \Vert Ax - b \Vert_M = \frac{1}{2}(b-Ax)^T M (b-Ax),
\end{eqnarray}
see for example~\cite{CE02}, Theorem 7.6. It follows that SIRT is both boundedly convergent and non-expansive, and therefore so is the operator $\mathbf{B}_w$ \eref{E:OS-SIRT2}, since it represents SIRT applied to a submatrix consisting of rows of $A$ and $b$. Since BI-SART~\eref{E:OS-SIRT1} is a composition of these operators together with projection onto the non-negative orthant, which also satisfies both properties, BI-SART satisfies both properties as well. Our proximity function $Pr(x) = \Vert Ax-b\Vert$ is also uniformly continuous (cf.~\cite{hgdc12}), and so therefore BI-SART is strongly perturbation resilient, guaranteeing convergence to a $\varepsilon$-compatible solution.

\section{Numerical Experiments}\label{S:exp}
\subsection{Experimental setup}
To validate the PnP-sup approach, we performed two sets of numerical experiments, one simulating sparse-view CT, and the other simulating low-dose CT. Both sets of experiments used $512 \times 512$ pixel CT image slices from a lung CT study obtained from the Cancer Imaging Archive~\cite{TCIA,QIN_LUNG_CT}, which included images from forty-two patients ranging from the lower abdomen to upper thorax. The ASTRA toolbox~\cite{PBS11,VPCJB16} was used both to generate simulated fan-beam sinogram data from the CT images, and also to implement the basic algorithm \eref{E:OS-SIRT2} and its superiorized variants in Python. The baseline normal-dose CT sinogram for each image consisted of 900 equally spaced views over 360$^\circ$, with a beam intensity of $I_0 = 10^6$ photons per view. The pixel size was assumed to be 0.0568 mm.

In both sets of experiments, we compared the results of the PnP-sup algorithm with both the basic algorithm, and two conventionally superiorized algorithms using total variation (TV) as the penalty function $\phi$. We define the TV function as
\begin{eqnarray}
\phi(x) = \sum_{m,n} \sqrt{ (x_{m+1,n}-x_{m,n})^2 + (x_{m,n+1}-x_{m,n})^2 + \epsilon^2},\label{E:TV}
\end{eqnarray}
where $x$ is a doubly-indexed 2D image, and $\epsilon=10^{-6}$ is a small parameter to ensure that $\phi$ is differentiable everywhere, so that the negative gradient can be used as a non-ascending direction in Algorithm~\ref{A:sup_conv}. Total variation is a frequently used penalty function in low-dose and sparse-view CT imaging, and within the superiorization literature, including \cite{hgdc12} and some of our previous work \cite{hwf17,humphries2020comparison}. The two TV-superiorized methods were denoted as follows:
\begin{itemize}
    \item BI-SART-TV: This approach was implemented following \aref{A:sup_conv}, using values of $\gamma = 0.9995$ and $N=20$, as these were found to give good results empirically. These parameter values are similar to those used in other studies on CT imaging using total variation, e.g.~\cite{hgdc12,humphries2020comparison,hwf17}.
    \item BI-SART-TVa: This approach uses the adaptive step size proposed in~\cite{naae22}, as given in \aref{A:sup_adap}. Following the paper, we used parameters $\displaystyle\alpha_0 = \frac{\phi(\tilde{x})}{2}$, $\displaystyle \epsilon = \frac{\phi(\tilde{x})}{200}$, where $\tilde{x}$ is a reconstructed image generated by one iteration of BI-SART, beginning from $x^0 = 0$.
\end{itemize}

For both the PnP-sup and TV-superiorized methods, the proximity function was defined to be the 2-norm residual:
\begin{equation}
    Pr_T(x) = \norm{A x - b}.\label{E:prox}
\end{equation}

In addition to a visual/qualitative assessment, we also judged the quality of reconstructed images using three quantitative metrics: peak signal-to-noise ratio (PSNR), Structural SIMilarity index (SSIM)\cite{SSIM}, and relative error in total variation ($\Delta TV \%$). The PSNR of a reconstructed image $x$ versus the baseline image $y$ is computed as
\begin{eqnarray}
PSNR(x,y) &= 10 \log_{10} \left( \frac{y_{max}}{MSE(x,y)} \right), \textrm{where} \label{E:PSNR} \\MSE(x,y) &= \frac{1}{N}\sum_j \left(x_{j} - y_{j}\right)^2, \: \textrm{ and } y_{max} = \max_j \{ y_j \}. 
\end{eqnarray}
Thus, larger PSNR values indicate better agreement with the true image. The SSIM is calculated as
\begin{eqnarray}
SSIM(x,y) &= \frac{ (2 \mu_x \mu_y + C_1)(2 \sigma_{xy} + C_2)}{(\mu_x^2 + \mu_y^2 + C_1)(\sigma_x^2 + \sigma_y^2 + C_2)},\label{E:SSIM}
\end{eqnarray}
where $\mu_x, \mu_y$ are the means of the two images, and $\sigma_x, \sigma_y$ and $\sigma_{xy}$ the variances and covariance, respectively. The small constants $C_1$ and $C_2$ are included to avoid instability when the denominator is near zero. The SSIM produces a score between 0 and 1, with 1 corresponding to identical images. It is argued that SSIM provides a better score of perceptual similarity between images than measures based on mean square error, such as PSNR~\cite{SSIM}; both are frequently used to assess image quality. The relative error in total variation was defined as
\begin{eqnarray}
    \Delta TV \%(x,y) = \frac{\phi(y)-\phi(x)}{\phi(y)} \times 100\%,
\end{eqnarray}
with $\phi$ given in \eref{E:TV}. A negative value therefore indicates that the TV of the reconstructed image was larger than the true image, and vice-versa.

\subsection{Plug-and-play superiorization with a deep neural network}\label{S:exp:nn}

For this experiment, we considered the use of a convolutional neural network (CNN) to improve the quality of images generated from sparse-view CT data using PnP-sup. The use of CNNs for sparse-view CT, and many other challenging scenarios in CT imaging, has attracted considerable research interest in the last ten years. There are a number of ways in which the power of neural networks can be leveraged for CT image reconstruction. Some examples include:
\begin{enumerate}
    \item Post-processing, in which an image is first reconstructed using a conventional algorithm, such as filtered backprojection or SART, and a network is trained to remove the resulting artifacts~\cite{zldxc18,hy18,hscsx19}.
    \item ``Unrolled'' algorithms, in which some predetermined number of iterations of an algorithm, such as SART or gradient descent, are unrolled as a CNN that includes trainable layers as well~\cite{czczz18,AO17,JMMTSH22}.
    \item Plug-and-play type algorithms, where a neural network is trained separately from an iterative algorithm such as gradient descent or ADMM, then inserted into the iteration~\cite{clpk17,gupta2018cnn}.
    \item Full inversion methods, which train a network to reconstruct an image directly from a sinogram~\cite{LLZMC19}.
\end{enumerate}
Using a neural network within the PnP-sup framework falls under the third category, where the network is trained apart from the iterative algorithm, then incorporated into the iteration. For our experiments, we used a CNN architecture previously proposed in~\cite{zzcmz17} for denoising of photographic images. The network consists of 17 convolutional layers, each consisting of sixty-four $3\times 3$-pixel convolutional filters, batch normalization, and rectified linear unit (ReLU) activation. It attempts to learn a residual mapping $\Psi(x)$ such that $\Psi(x) \approx y - x$, where $y$ is the ground truth and $x$ is the noisy or otherwise degraded image; in other words, rather than learning a mapping from a noisy image to a clean image, the network maps the noisy image to the actual pattern of noise (which can then be subtracted to obtain the clean image). This network was studied in our previous work~\cite{hscsx19} as a post-processing operator for sparse-view and low-dose CT. Code for the network is available at \url{https://github.com/cszn/DnCNN}, and a more thorough description of the architecture can be found in~\cite{zzcmz17}.

\begin{figure}
    \centering
    \includegraphics[width=\linewidth]{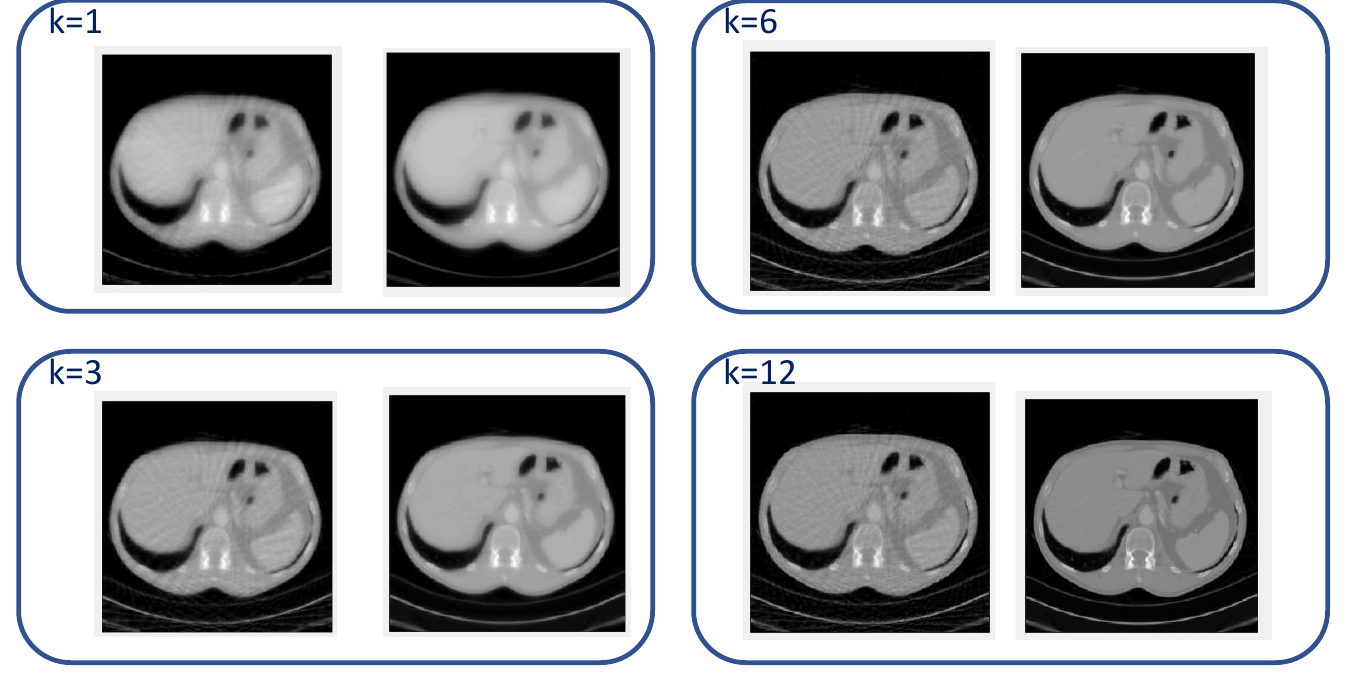}
    \caption{Example of paired sparse-view (left images) and normal-dose (right images) training data at different stages of iterative reconstruction ($k=1,3,6, 12$ iterations). }
    \label{F:trainingdata}
\end{figure}

To train this network for use within PnP-sup, we generated training data by simulating normal-dose (900 views) fan-beam sinograms from 1000 CT image slices, together with sparse-view (60 views) sinograms from the same image slices. We chose 60 views for the sparse-view data as this is comparable to other papers in the literature, e.g.~\cite{zldxc18,czczz18,JMMTSH22}. Both sinograms were then used to reconstruct paired images using 12 iterations of BI-SART~\eref{E:basic} with 10 subsets. Because the network is applied after every iteration of the basic algorithm, we used image pairs corresponding to the first, third, sixth and twelfth iterations as training data. This ensured that the network was trained on data corresponding both to early and later iterations of the algorithm, since the characteristics of the image vary across iterations, as illustrated in~\Fref{F:trainingdata}. Thus there were a total of 4000 full-size CT images used to train the network; however, the actual number of training samples was much larger, as the network is trained on small $32 \times 32$ pixel patches extracted from these images~\cite{zzcmz17}. We refer to this algorithm as PnP-sup-NN.

We used a set of 20 CT image slices from a separate patient as the test data set to assess the performance of the method. First, an image was reconstructed from a 60-view sinogram using the basic algorithm, defined to be 12 iterations of BI-SART with 10 subsets, just as for the training data. Evaluating the proximity function $Pr_T$ for this image provided the value of $\varepsilon$ used as the stopping criterion for BI-SART-TV, BI-SART-TVa, and PnP-sup-NN. For PnP-sup-NN, a value of $\gamma = 0.95$ was used (cf. Algorithm~\ref{A:sup_pnp}); this value was smaller than that used by BI-SART-TV, to ensure that the CNN does not introduce large perturbations at later iterations of the algorithm. As a fourth point of comparison, we took the same neural network architecture and trained it as a post-processing operation by only using normal-dose and sparse-view paired data corresponding to 12 iterations of BI-SART (images resembling those in the bottom right of \Fref{F:trainingdata}). This is similar to the approach used in our earlier work~\cite{hscsx19}.

\subsection{Plug-and-play superiorization with a black-box denoiser}

In the second set of experiments, we applied perturbations generated by the black-box denoiser BM3D \cite{dabov2007image,MAF20} within the PnP-sup framework to improve the quality of images generated from low-dose CT data. We refer to this algorithm as PnP-sup-BM3D. Since there was no need for training data in this case, the same 20 CT image slices used for test data in the previous set of experiments were used to evaluate the approach here. In addition to the normal-dose data ($I_0 = 1\ee6$), three reduced intensity levels were simulated: $I_0 = 5\ee{4}$, $I_0 = 2.5\ee{4}$, and $I_0 = 1\ee{4}$. To generate values of $\varepsilon$ as the stopping criterion for both BI-SART-TV and PnP-sup-BM3D, the BI-SART algorithm was run for 18, 12, or 8 iterations, corresponding, respectively, to the three noise levels listed previously. The number of iterations was reduced for higher noise levels due to the well-known semiconvergence property~\cite{EHN12}, which causes image quality to deterioriate in the presence of noisy data, if the algorithm is not terminated early; in general, the noisier the data, the earlier the algorithm should be stopped.

BI-SART-TV and BI-SART-TVa were run with the same parameters as described in the previous section. For PnP-sup-BM3D, we found that we needed to modify the iteration given in Algorithm~\ref{A:sup_pnp} to obtain good results; in particular, two additional parameters were introduced:
\begin{itemize}
    \item $k_{min} \geq 0$: an integer indicating when to begin applying perturbations to the basic algorithm, and
    \item $k_{step} \geq 1$: an integer indicating how many iterations of the basic algorithm to run in between each perturbation.
\end{itemize}
The first to fourth lines of the loop (from $z \gets \Psi(x^k)$, to $\displaystyle \beta_k \gets \min \{\alpha \gamma^\ell, \norm{v^k} \}$ ) were then executed if and only if the following condition was satisfied:
\begin{equation}
(k \geq k_{min}) \textrm{ and } ((k - k_{min})~\%~k_{step} == 0) \label{E:kstep}
\end{equation}
The reason for this modification was that we found that applying BM3D to initial iterates of the solution was not useful, as noise had not yet begun to emerge in the reconstructed image (as expected, due to the semiconvergence property). Furthermore, if BM3D was applied in between every iteration of the basic algorithm, the perturbations introduced were initially so large that it took many iterations before the residual, $Pr_T(x^{k+1})$, became small enough for the algorithm to terminate; unfortunately, by this point, the perturbations introduced by BM3D has been damped to the extent that the final iterates produced still contained significant noise, and were barely any better than those produced by the basic algorithm. Introducing the parameters $k_{min}$ and $k_{step}$ alleviated these respective issues. 

A second modification was also made, to initialize the starting step size, $\alpha$, to be the size of the first perturbation, $\alpha = \Vert v^{k_{min}} \Vert$, rather than making it a user-specified parameter. This allowed a ``full step'' of BM3D to be applied when first denoising the image, rather than potentially damping it if $\alpha$ is chosen to be too small. Values of $k_{min} = 15, 10,$ and $5$ were used for the three test cases (in order from least to most noisy data), with $k_{step} = 5$ in the first two cases, and $4$ in the last (most noisy) case; $\gamma = 0.75$ was used as the rate of damping for the perturbations. Finally, BM3D was also applied as a post-processing step to the iterations produced by the basic algorithm (BI-SART), as a third point of comparison.

\section{Results}\label{S:results}

\subsection{Plug-and-play superiorization with a deep neural network}
\begin{table}
\caption{Results of experiments using PnP-sup-NN for sparse-view CT. Second column shows average and standard deviation of PSNR values over the 20 test images for each method; third column the corresponding values for SSIM, fourth column the average $\Delta TV \%$ values, fifth column the average number of iterations of the algorithm needed to generate the result, sixth column the average runtime needed for each method, and the last column the average value of $Pr_T(x) = \Vert Ax-b \Vert$ attained by the solution. Metrics showing best performance according to one-way ANOVA are bolded.}\label{T:exp1}
\footnotesize
\begin{tabular}{lcccccc}
\hline
Method &PSNR &SSIM &$\Delta TV\%$&Iterations &$t$ (s)&$\varepsilon$ \\
\hline
BI-SART 		&$33.96 \pm 2.53$ &$0.862 \pm 0.040$ &-28.9&12 &6&1.71 \\
BI-SART-TV 	&$\mathbf{37.54 \pm 2.45}$ &$\mathbf{0.950 \pm 0.014}$&20.8&68 &259&1.70 \\
BI-SART-TVa &$36.59 \pm 2.34$          &$\mathbf{0.935 \pm 0.019}$                           &$-10.0$   &61   &50   &1.70\\
PnP-sup-NN 	&$\mathbf{38.89 \pm 2.02}$ &$\mathbf{0.949 \pm 0.013}$&-0.2&14 &9&1.67 \\
NN-Post	&$\mathbf{37.97 \pm 2.09}$ &$0.913 \pm 0.024$&6.5&12 &6&5.78 \\
\hline
\end{tabular}
\end{table}

Quantitative results for the sparse-view CT experiments using PnP-sup-NN are shown in Table~\ref{T:exp1}. We can see that all of BI-SART-TV, BI-SART-TVa, PnP-sup-NN, and post-processing by the neural network (NN-Post) substantially improve the quality of the image versus the basic algorithm (BI-SART), as measured both by PSNR and SSIM. A one-way ANOVA test indicated that BI-SART-TV and PnP-sup-NN produced comparable results to each other with respect to both PSNR and SSIM (no statistically significant difference), while BI-SART-TVa provided comparable SSIM values only to these top two methods, and NN-Post comparable PSNR values only. The $\Delta TV \%$ values indicate that PnP-sup-NN provided the best results, with a discrepancy of only -0.2\%; The values of 20.8\% and 6.5\% for BI-SART-TV and NN-Post, respectively, indicates that those methods produced images with lower TV than the true image, most likely due to oversmoothing. BI-SART-TVa, on the other hand, did not have as strong a smoothing effect, as indicated by the negative $\Delta TV\%$ value.

One noticeable difference between PnP-sup-NN and the two TV-superiorized methods is that the PnP-sup-NN method requires far fewer iterations to converge to an $\varepsilon$-compatible solution. This resulted in significant computational savings, as indicated by the runtime values in Table~\ref{T:exp1}; especially when compared to BI-SART-TV. (Iterations of BI-SART-TVa are less expensive compared to iterations of BI-SART-TV, due to the elimination of the inner loop over $N$). Intuitively, this makes sense; the neural network is trained to map sparse-view reconstructed images at various points of the iterative process to their corresponding normal-dose (baseline) images, which are also obtained by attempting to fit the measured data. Thus in some sense, the basic algorithm and the superiorization step do not compete with each other. In contrast, when superiorizing with respect to TV, there is an inherent conflict between trying to fit the data and trying to reduce the TV of the image, which results in a longer time to convergence. 

Additionally, we can see that while both the TV-superiorized methods and PnP-sup-NN are (as expected) able to provide solutions meeting the same level of constraints compatibility ($\varepsilon$) as the BI-SART algorithm, the images produced by NN-Post have much higher discrepancies with the sinogram data. This is not surprising given that the neural network is applied as a post-processing operation in this case, which has no knowledge of the original measured data.

\begin{figure}
    \includegraphics[width=\linewidth]{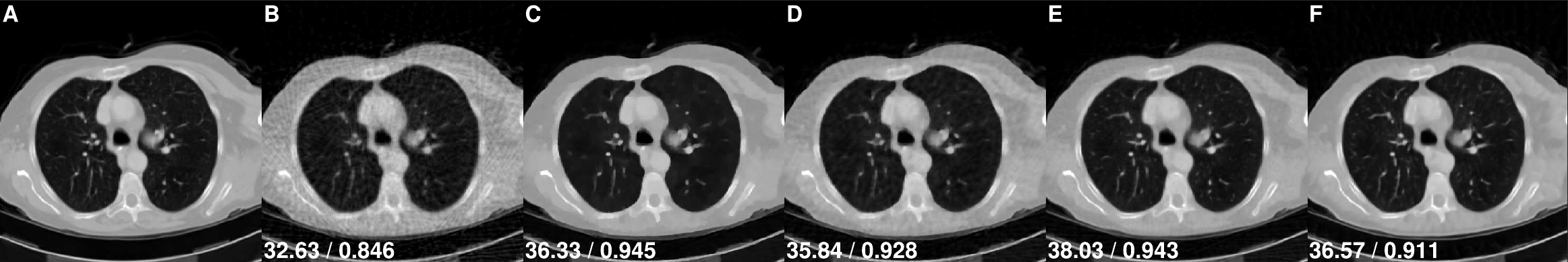}
    \includegraphics[width=\linewidth]{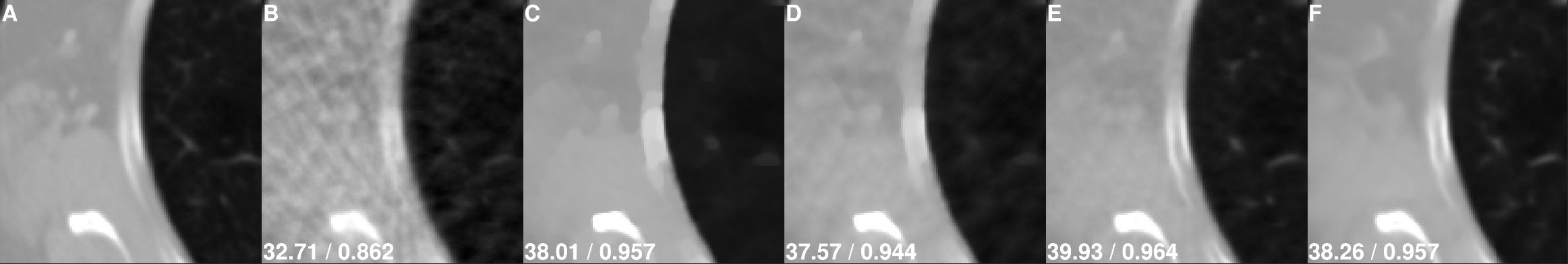}
    \includegraphics[width=\linewidth]{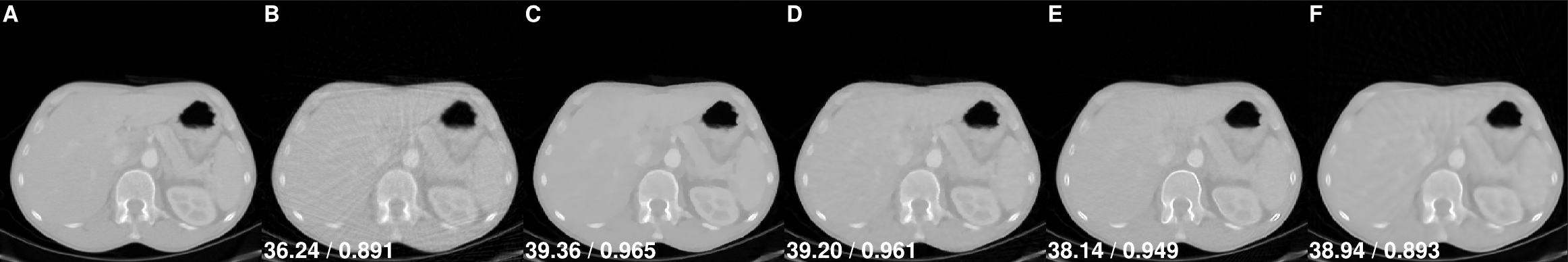}
    \includegraphics[width=\linewidth]{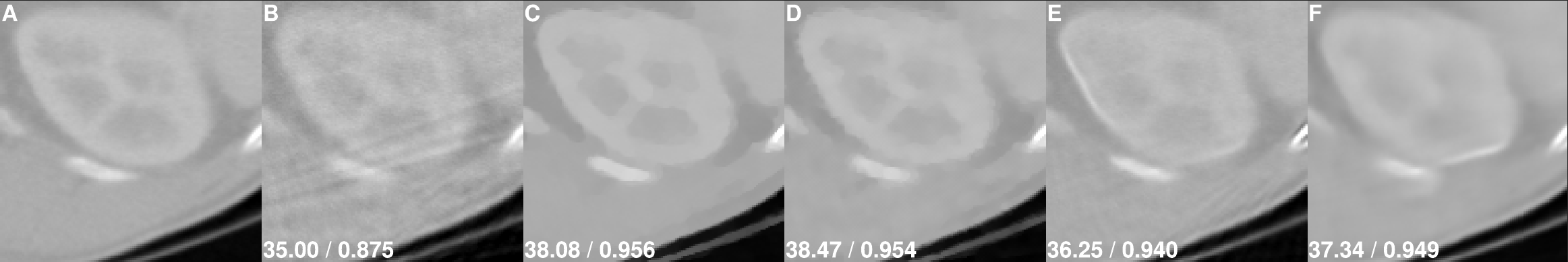}
    \caption{Representative sparse-view images reconstructed using different approaches. Column A: ground truth image; Column B: image reconstructed using BI-SART; Column C: image reconstructed using BI-SART-TV; Column D: image reconstructed using BI-SART-TVa; Column E: image reconstructed using PnP-sup-NN; Column F: image post-processed by neural network. Second and fourth rows show zoomed-in region indicated in red box on ground truth image above the row. Color window is [0, 0.3] cm$^{-1}$. PSNR/SSIM values for each image or zoomed-in region shown in bottom left.}\label{F:PnPS-NN images}
\end{figure}

Figure~\ref{F:PnPS-NN images} shows two representative images from the test data set, the first of the patient's lungs, the second of the abdomen. In both cases the TV-superiorized images (columns C and D) shows the characteristic ``cartoonish'' texture sometimes associated with the use of TV, and it is clear that some structures such as the rib bones have been over-smoothed. The texture of the images reconstructed using PnP-sup-NN, on the other hand, more closely resembles that of the ground truth image. That being said, it is apparent that some false artifacts have been introduced by the neural networks in the abdominal CT image (bottom row, columns E and F), as some edges of the organ in the highlighted region have been enhanced unrealistically. These artifacts could potentially be alleviated by better training of the network. This is also reflected in the PSNR/SSIM values shown in the bottom left of each image; for the lung image, PnP-sup-NN provides the best results for both the overall image and the zoomed-in region, while for the abdominal image, BI-SART-TV and BI-SART-TVa gave better results.

\subsection{Plug-and-play superiorization with a black-box denoiser}

\begin{table}
\caption{Results of experiments using PnP-sup-BM3D for low-dose CT at the three different noise levels. First column indicates the noise level for the experiment; second column shows average and standard deviation of PSNR values over the 20 test images for each method; third column the corresponding values for SSIM, fourth column the average $\Delta TV \%$ values, fifth column the average number of iterations of the algorithm needed to generate the result, sixth column the average runtime needed for each method, and the last column the average value of $Pr_T(x) = \Vert Ax-b \Vert$ attained by the solution. Metrics showing best performance according to one-way ANOVA are bolded.}\label{T:exp2}
\footnotesize
\begin{tabular}{llcccccc}
\hline
Exp. & Method &PSNR &SSIM&$\Delta TV \%$&Iterations &$t$ (s)&$\varepsilon$ \\
\hline
{\tt 5e4}&BI-SART 		&$38.98 \pm 1.92$ &$0.923 \pm 0.037$ &$-48.5$&18 &1 &16.17  \\
&BI-SART-TV 	&$\mathbf{40.08 \pm 2.52}$ &$\mathbf{0.950 \pm 0.040}$ &$-14.9$&337 &394 &16.16 \\
&BI-SART-TV-a   &$\mathbf{39.86 \pm 2.62}$ &$\mathbf{0.949 \pm 0.040}$ &$-14.1$ &677 &128 &16.17\\
&PnP-sup-BM3D 	&$\mathbf{41.00 \pm 2.21}$ &$\mathbf{0.944 \pm 0.038}$&$-35.3$&61 &33&16.13 \\
&BM3D-Post	&$38.65 \pm 2.20$ &$\mathbf{0.960 \pm 0.012}$&$38.5$&18 &5 &19.68\\
\hline
{\tt 2.5e4}&BI-SART 		&$36.25 \pm 2.02$ &$0.928 \pm 0.025$ &-12.0&12 &1 &26.04 \\
&BI-SART-TV 	&$37.89 \pm 2.64$ &$\mathbf{0.950 \pm 0.021}$ &36.3&71 &107 &26.02 \\
&BI-SART-TV-a   &$38.04 \pm 2.63$ &$\mathbf{0.949 \pm 0.022}$ &6.5 &51 &10 &25.96\\
&PnP-sup-BM3D 	&$38.01 \pm 2.65$ &$\mathbf{0.954 \pm 0.020}$&18.2&18 &7& 25.91\\
&BM3D-Post	&$36.35 \pm 2.24$ &$\mathbf{0.946 \pm 0.016}$&43.7&12 &4 &28.01\\
\hline
{\tt 1e4}&BI-SART 		&$33.91 \pm 2.06$ &$0.900 \pm 0.030$&-11.6 &8 &1 &41.87 \\
&BI-SART-TV 	&$35.27 \pm 2.67$ &$\mathbf{0.929 \pm 0.024}$&40.1&37 &58 &41.75 \\
&BI-SART-TV-a   &$35.31 \pm 2.41$ &$\mathbf{0.930 \pm 0.023}$ &12.0 &18 &4 &41.60\\
&PnP-sup-BM3D 	&$35.25 \pm 2.28$ &$\mathbf{0.932 \pm 0.022}$&20.9&11 &7 &41.42 \\
&BM3D-Post	&$34.07 \pm 2.25$ &$\mathbf{0.927 \pm 0.021}$&49.0&8 &4&43.87 \\
\hline
\end{tabular}
\end{table}
Quantitative results for the low-dose CT experiments are shown in Table~\ref{T:exp2}. We can observe that although the PSNR values were typically improved somewhat by applying BI-SART-TV, BI-SART-TVa, PnP-sup-BM3D, or BM3D-Post versus BI-SART, in only one instance was the improvement found to be statistically significant, for the least noisy dataset. For the SSIM values, all four methods tended to provide comparable improvement. The positive $\Delta TV\%$ values for the two higher-noise experiments indicate that all of the approaches tended to produce images with substantially lower TV values than the true image, indicating possible oversmoothing; of these, SART-TVa and PnP-sup-BM3D resulted in values substantially closer to the true image than the other two approaches. In the low-noise experiment, most of the $\Delta TV \%$ values were negative, on the other hand, with BI-SART-TV and BI-SART-TVa having the smallest $\Delta TV \%$ values.

As with the sparse-view experiments, when comparing BI-SART-TV and BI-SART-TVa to PnP-sup-BM3D, there is a notable difference in the total number of iterations required to converge, with the former two methods requiring significantly more iterations, especially for the two less noisy datasets. We note that unlike the previous set of experiments, where perturbing iterates using the neural network only slightly increased the total number of iterations required, doing so with BM3D does tend to require more iterations to reach the same level of constraints-compatibility as BI-SART. This indicates that the denoising applied by BM3D does compete with the goal of constraints compatibility, to some extent. Applying TV superiorization, however, tended to increase the required number of iterations by a much greater factor. This is likely due to the fact that the perturbations introduced by PnP-sup-BM3D were applied less frequently, due to the introduction of the parameters $k_{min}$ and $k_{step}$, and decayed in size more quickly, due to the smaller value of $\gamma$ used. The runtime values in Table~\ref{T:exp2} show that BI-SART-TV was by far the most computationally intensive algorithm, while BI-SART-TVa required substantially lower runtime; even less than PnP-sup-BM3D in the experiment with the noisiest dataset.

Finally, in comparing the $\varepsilon$ (residual) values in the last column of the table, we can see again that all of BI-SART-TV, BI-SART-TVa and PnP-sup-BM3D eventually attain solutions with an equivalent, or smaller value than BI-SART. When the image produced by BI-SART is postprocessed using BM3D (BM3D-Post), the residual value is greater than that of the image produced BI-SART, though not to the same extent as in the previous set of experiments that used the neural network. In comparison to Table~\ref{T:exp1}, we note that the residual values are much greater here due to the larger size of the projection data (900 views vs 60 views). The value of $\varepsilon$ at termination also increases as the noise level increases, due both to the lower number of iterations of BI-SART that were run as noise increased, as well as the greater inconsistency introduced to the data by higher noise levels.

\begin{figure}
    \includegraphics[width=\linewidth]{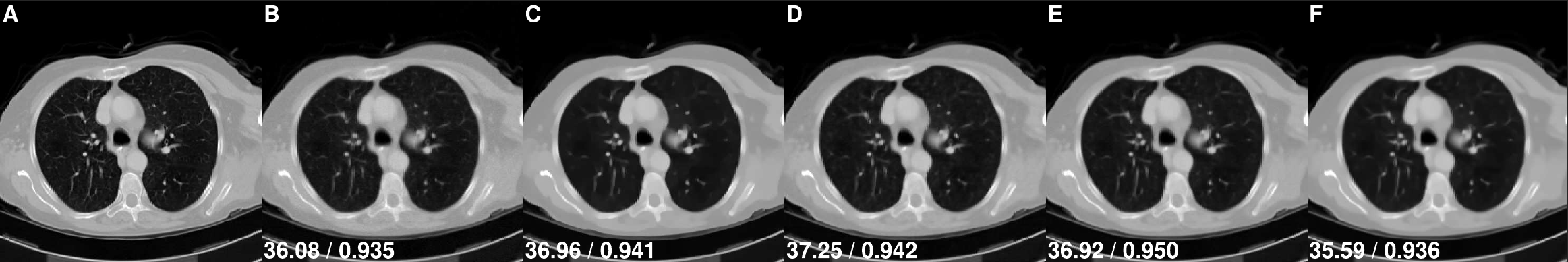}
    \includegraphics[width=\linewidth]{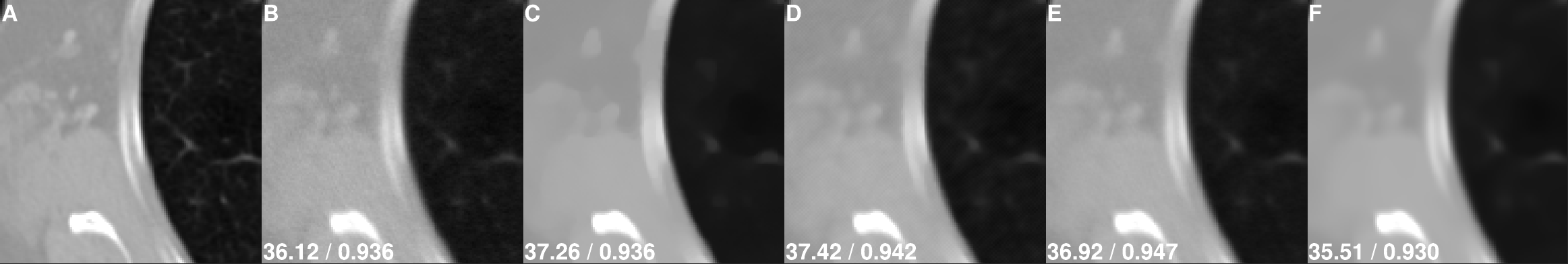}
    \includegraphics[width=\linewidth]{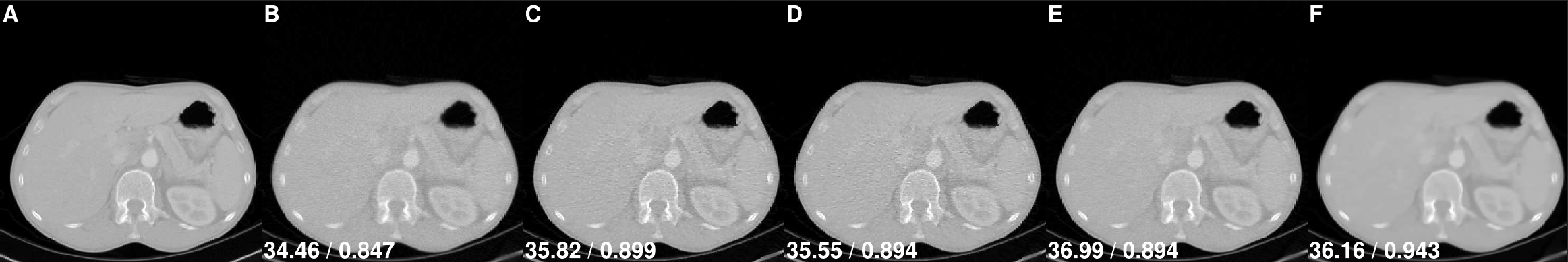}
    \includegraphics[width=\linewidth]{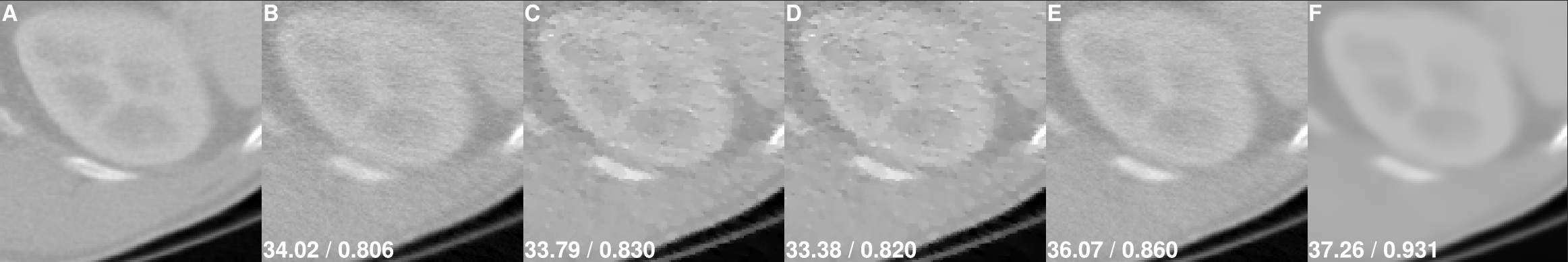}
    \caption{Representative images at dose level of $I_0 = 2.5 \times 10^4$ reconstructed using different approaches. Column A: ground truth image; Column B: image reconstructed using BI-SART; Column C: image reconstructed using BI-SART-TV; Column D: image reconstructed using BI-SART-TVa; Column E: image reconstructed using PnP-sup-BM3D; Column F: image post-processed by BM3D. Second and fourth rows show zoomed-in region indicated in red box on ground truth image above the row. Color window is [0, 0.3] cm$^{-1}$. PSNR/SSIM values for each image or zoomed-in region shown in bottom left.}\label{F:PnPS-BM3D images}
\end{figure}

Figure~\ref{F:PnPS-BM3D images} shows the same two representative images as Figure~\ref{F:PnPS-NN images} for the simulation at a dose level of $I_0 = 2.5 \times 10^4$. We first note that while the dose level is the same for both images, the noise is much more apparent in the abdominal image; this is because the noise is Poisson-distributed proportionally to the {\em measured} intensity on the other side of the object, which is lower for the abdominal image, as it consists of more attenuating material than the image of the lungs. With respect to the performance of the different reconstruction approaches, we observe that BI-SART-TV and BI-SART-TVa give inconsistent results; the images of the lungs have been over-smoothed, while the images of the abdomen appears even noisier than those produced by BI-SART. The images produced by BM3D-Post, on the other hand, both appear oversmoothed. Finally, the lung image produced by PnP-sup-BM3D has provided denoising of the image without oversmoothing, while for the abdominal image, there is a less noticeable improvement. The PSNR and SSIM values for the lung image show that PnP-sup-BM3D provides comparable values to the images produced by BI-SART-TV and BI-SART-TVa. For the abdominal image, it provides the highest PSNR value for the overall image, but lower PSNR than BM3D-Post for the zoomed-in portion; BM3D-post provides the highest SSIM in both cases, despite the overly smooth appearance of the image.

\section{Discussion}\label{S:disc}
We can make several observations based on the results of the numerical experiments shown in the previous section. The first is that in both sets of experiments, the PnP-sup approach is generally successful in improving the quality of the solutions found by the basic algorithm. The average PSNR and SSIM values of the test data set were higher in every instance, although ANOVA testing indicated that the differences in PSNR were not statistically significant in two of the low-dose scenarios. The appearance of the images is also generally better, especially in the sparse-view experiments. We note that while the improvements in image quality provided by PnP-sup are sometimes not dramatic, we have not focused on optimizing performance of the PnP-sup approach in this paper, but rather view it as providing proof-of-concept for the methodology. It is likely that improved results could be obtained, for example, by using a more sophisticated neural network architecture in Section~\ref{S:exp:nn}, or tuning other algorithmic parameters. At the very least, we can say that PnP-sup provides improved or comparable image quality to all of the approaches we compared against, while providing significant computational savings over the BI-SART-TV method, and better data fidelity than the post-processing approaches. The BI-SART-TVa method provides the closest results, but still underperforms PnP-sup in the sparse-view experiments, and was also computationally more expensive in two of the four scenarios we considered.

Comparing PnP-sup with the conventionally superiorized BI-SART-TV and BI-SART-TVa, we see that the average PSNR and SSIM values obtained using the plug-and-play approach were comparable to or higher than those obtained using the two TV-superiorized methods in every experiment, though they were often quite close. The images produced using PnP-sup tended to have better appearance, and to avoid the oversmooth or blocky texture that images produced by TV minimization can exhibit. Additionally, the $\Delta TV\%$ we computed indicated that the PnP-sup approaches we used were effective in reducing the average TV values of the reconstructed images when compared to BI-SART, even though TV was not explicitly minimized by the method. In the sparse-view experiments, the images it produced also more closely matched the TV of the true images than the images produced by BI-SART-TV and BI-SART-TVa.

The PnP-sup method used in our two sets of experiments offered significant computational savings over BI-SART-TV, both in terms of the number of iterations required for convergence, and total runtime. This is likely due in part to a somewhat conservative choice of parameters for BI-SART-TV ($\gamma = 0.9995$ and $N=20$), meaning that 20 perturbations were introduced between each iteration of BI-SART, and that the perturbations decreased slowly in size. Smaller values of these two parameters would result in fewer iterations needed for convergence, but based on our experience, might also produce lower-quality images. The BI-SART-TVa method, which uses an adaptive stepsize, was computationally less expensive than BI-SART-TV while providing comparable results in most cases, and also provided comparable runtime to PnP-sup-BM3D in two of the lose-dose scenarios.

Although we have not conducted a detailed sensitivity analysis in this paper, it is likely that the PnP-sup method (\aref{A:sup_pnp}) suffers from some of the same shortcomings as the conventionally superiorized \aref{A:sup_conv}, namely, sensitivity to the choice of parameter $\gamma$ controlling stepsize, stopping criterion $\varepsilon$, and the parameters $k_{min}$ and $k_{step}$ introduced in \eref{E:kstep}, if they are used. These parameters may need to be fine-tuned for specific applications. It would be of interest to investigate approaches to determine these parameters automatically or adaptively, similar to the stepsize rule used in \aref{A:sup_adap}.



Post-processing the images produced by the basic algorithm using the neural network or BM3D also tended to improve image quality, and provide PSNR and SSIM values similar to, or slightly lower, than those obtained by the corresponding PnP-sup algorithms. As shown in the numerical experiments, the main advantage of the latter is the guarantee of $\varepsilon$-compatibility provided by the superiorization methodology (SM), which is not the case when a neural network or BM3D is applied as a postprocessing operation. This is important in contexts such as CT imaging where $\varepsilon$-compatibility reflects agreement with the data measured by the camera system.

In earlier work by one of the authors with colleagues~\cite{JMMTSH22}, we presented an unrolled neural network approach inspired by the SM. It is interesting to compare this unrolled approach with the plug-and-play approach presented here, especially in the context where the perturbations are generated by a neural network. A similarity between the two approaches is that, contrary to the typical use of superiorization, the perturbations are generated by a neural network, rather than as nonascending directions of some prescribed function $\phi(x)$. The sense in which the iterates are expected to be superior to those produced by the basic algorithm is therefore less straightforward to assess.  In both cases, however, if the network is trained appropriately, we do observe significant improvements in terms of standard image quality metrics such as PSNR and SSIM, even though these are not explicitly being reduced by the perturbations.

It is also interesting to compare the differences between the two approaches and their relative advantages and disadvantages. A main advantage of the plug-and-play approach is that it fits more closely with the existing SM, as the superiorized algorithm can be run indefinitely until converging to a $\varepsilon$-compatible solution. On the other hand, no such guarantee of $\varepsilon$-compatibility exists for the output of the unrolled algorithm, since the total number of iterations of the algorithm must be fixed in order to train the entire iterative process as a neural network. On the other hand, since the entire iterative process is included in the training of the unrolled algorithm (including the operations applied as part of the basic iteration), it may be able to better optimize performance versus the plug-and-play approach, in which the neural network is blind to the action of the basic algorithm. We note that~\cite{JMMTSH22} included a 60-view sparse-view experiment as well, in which the unrolled algorithm was able to achieve somewhat better results than those observed in Section~\ref{S:exp:nn} of this paper.

As discussed in the aforementioned paper, one possibility is to take the output of the unrolled network and use it as the initial iterate for the basic, conventionally superiorized, or even plug-and-play superiorized algorithm, which could be iterated (hopefully for only a small number of iterations) until a $\varepsilon$-compatible solution is reached. Even though this would be a two-step process, one could still view the entire process overall as a superiorization of the basic algorithm, consisting of finitely many perturbations of unknown size (generated by the unrolled portion), followed by a series of perturbations of decreasing size. Such a hybridization might provide the best of both worlds.     

\section{Conclusions}\label{S:concl}

In this paper we present a novel framework for plug-and-play superiorization, in which iterates produced by a basic algorithm are perturbed by a black-box operation (such as a denoiser, or neural network). This is in contrast to the conventional superiorization methodology (SM), in which the perturbations are made in nonascending directions of some specified function. By controlling the size of the perturbations appropriately, we are able to provide the same guarantees of convergence provided by the conventional SM. Our numerical experiments on two problems in CT imaging demonstrate the promise of the plug-and-play approach. We expect that this framework may find use in other applications as well, especially given the ability it provides to incorporate neural networks into existing algorithms.


\section*{References}

\begin{thebibliography}{41}
\expandafter\ifx\csname natexlab\endcsname\relax\def\natexlab#1{#1}\fi
\providecommand{\url}[1]{\texttt{#1}}
\providecommand{\href}[2]{#2}
\providecommand{\path}[1]{#1}
\providecommand{\DOIprefix}{doi:}
\providecommand{\ArXivprefix}{arXiv:}
\providecommand{\URLprefix}{URL: }
\providecommand{\Pubmedprefix}{pmid:}
\providecommand{\doi}[1]{\href{http://dx.doi.org/#1}{\path{#1}}}
\providecommand{\Pubmed}[1]{\href{pmid:#1}{\path{#1}}}
\providecommand{\bibinfo}[2]{#2}
\ifx\xfnm\relax \def\xfnm[#1]{\unskip,\space#1}\fi
\bibitem[{Censor et~al.(2014)Censor, Davidi, Herman, Schulte, and
  Tetruashvili}]{cdhst14}
\bibinfo{author}{Y.~Censor}, \bibinfo{author}{R.~Davidi},
  \bibinfo{author}{G.~T. Herman}, \bibinfo{author}{R.~W. Schulte},
  \bibinfo{author}{L.~Tetruashvili},
\newblock \bibinfo{title}{Projected subgradient minimization versus
  superiorization},
\newblock \bibinfo{journal}{J Optim. Theory Appl.} \bibinfo{volume}{160}
  (\bibinfo{year}{2014}) \bibinfo{pages}{730--747}.
\bibitem[{Guenter et~al.(2022)Guenter, Collins, Ogilvy, Hare, and
  Jirasek}]{gcohj22}
\bibinfo{author}{M.~Guenter}, \bibinfo{author}{S.~Collins},
  \bibinfo{author}{A.~Ogilvy}, \bibinfo{author}{W.~Hare},
  \bibinfo{author}{A.~Jirasek},
\newblock \bibinfo{title}{Superiorization versus regularization: A comparison
  of algorithms for solving image reconstruction problems with applications in
  computed tomography},
\newblock \bibinfo{journal}{Med Phys.} \bibinfo{volume}{49}
  (\bibinfo{year}{2022}) \bibinfo{pages}{1065--1082}.
\bibitem[{Davidi et~al.(2015)Davidi, Censor, Schulte, Geneser, and
  Xing}]{dcsgx15}
\bibinfo{author}{R.~Davidi}, \bibinfo{author}{Y.~Censor},
  \bibinfo{author}{R.~W. Schulte}, \bibinfo{author}{S.~Geneser},
  \bibinfo{author}{L.~Xing},
\newblock \bibinfo{title}{Feasibility-seeking and and superiorization
  algorithms applied to inverse treatment planning in radiation therapy},
\newblock \bibinfo{journal}{Contemp Math.} \bibinfo{volume}{636}
  (\bibinfo{year}{2015}) \bibinfo{pages}{83--92}.
\bibitem[{Jin et~al.(2013)Jin, Censor, and Jiang}]{jcj13}
\bibinfo{author}{W.~Jin}, \bibinfo{author}{Y.~Censor},
  \bibinfo{author}{M.~Jiang},
\newblock \bibinfo{title}{A heuristic superiorization-like approach to
  bioluminescence},
\newblock \bibinfo{journal}{International Federation for Medical and Biological
  Engineering (IFMBE) Proceedings} \bibinfo{volume}{39} (\bibinfo{year}{2013})
  \bibinfo{pages}{1026--1029}.
\bibitem[{Censor(2017)}]{c17}
\bibinfo{author}{Y.~Censor},
\newblock \bibinfo{title}{Can linear superiorization be useful for linear
  optimization problems?},
\newblock \bibinfo{journal}{Inverse Problems} \bibinfo{volume}{33}
  (\bibinfo{year}{2017}) \bibinfo{pages}{044006}.
\bibitem[{Censor(2024)}]{CensorSupPage}
\bibinfo{author}{Y.~Censor},
\newblock \bibinfo{title}{Superiorization and perturbation resilience of
  algorithms: A bibliography compiled and continuously updated}
  (\bibinfo{year}{2024}). \URLprefix
  \url{http://math.haifa.ac.il/YAIR/bib-superiorization-censor.html}.
\bibitem[{Censor et~al.(2017)Censor, Herman, and Jiang}]{CHJ17}
\bibinfo{author}{Y.~Censor}, \bibinfo{author}{G.~T. Herman},
  \bibinfo{author}{M.~Jiang},
\newblock \bibinfo{title}{Superiorization: theory and applications},
\newblock \bibinfo{journal}{Inverse Problems} \bibinfo{volume}{33}
  (\bibinfo{year}{2017}) \bibinfo{pages}{040301}. \URLprefix
  \url{https://doi.org/10.1088/1361-6420/aa5deb}.
  \DOIprefix\doi{10.1088/1361-6420/aa5deb}.
\bibitem[{Gibali et~al.(2020)Gibali, Herman, and Schn\"{o}rr}]{GHS20}
\bibinfo{author}{A.~Gibali}, \bibinfo{author}{G.~T. Herman},
  \bibinfo{author}{C.~Schn\"{o}rr},
\newblock \bibinfo{title}{Editorial: A special issue focused on superiorization
  versus constrained optimization: analysis and applications},
\newblock \bibinfo{journal}{Journal of Applied and Numerical Optimization}
  \bibinfo{volume}{2} (\bibinfo{year}{2020}) \bibinfo{pages}{1--2}. \URLprefix
  \url{https://doi.org/10.23952/jano.2.2020.1.01}.
  \DOIprefix\doi{10.1088/1361-6420/aa5deb}.
\bibitem[{Venkatakrishnan et~al.(2013)Venkatakrishnan, Bouman, and
  Wohlberg}]{vbw13}
\bibinfo{author}{S.~V. Venkatakrishnan}, \bibinfo{author}{C.~A. Bouman},
  \bibinfo{author}{B.~Wohlberg},
\newblock \bibinfo{title}{Plug-and-play priors for model based reconstruction},
\newblock in: \bibinfo{booktitle}{2013 IEEE Global Conference on Signal and
  Information Processing}, \bibinfo{organization}{IEEE}, \bibinfo{year}{2013},
  pp. \bibinfo{pages}{945--948}.
\bibitem[{Chan et~al.(2017)Chan, Wang, and Elgendy}]{cwe17}
\bibinfo{author}{S.~H. Chan}, \bibinfo{author}{X.~Wang}, \bibinfo{author}{O.~A.
  Elgendy},
\newblock \bibinfo{title}{Plug-and-play {ADMM} for image restoration:
  Fixed-point convergence and applications},
\newblock \bibinfo{journal}{IEEE Trans. Comput. Imag.} \bibinfo{volume}{3}
  (\bibinfo{year}{2017}) \bibinfo{pages}{84–98}.
\bibitem[{He et~al.(2018)He, Yang, Wang, Zeng, Bian, Zhang, Sun, Xu, and
  Ma}]{hywzb18}
\bibinfo{author}{J.~He}, \bibinfo{author}{Y.~Yang}, \bibinfo{author}{Y.~Wang},
  \bibinfo{author}{D.~Zeng}, \bibinfo{author}{Z.~Bian},
  \bibinfo{author}{H.~Zhang}, \bibinfo{author}{J.~Sun},
  \bibinfo{author}{Z.~Xu}, \bibinfo{author}{J.~Ma},
\newblock \bibinfo{title}{Optimizing a parameterized plug-and-play {ADMM} for
  iterative low-dose {CT} reconstruction},
\newblock \bibinfo{journal}{IEEE Trans. Med. Imag.} \bibinfo{volume}{38}
  (\bibinfo{year}{2018}) \bibinfo{pages}{371--382}.
\bibitem[{Censor et~al.(2010)Censor, Davidi, and Herman}]{cdh10}
\bibinfo{author}{Y.~Censor}, \bibinfo{author}{R.~Davidi},
  \bibinfo{author}{G.~T. Herman},
\newblock \bibinfo{title}{Perturbation resilience and superiorization of
  iterative algorithms},
\newblock \bibinfo{journal}{Inverse Problems} \bibinfo{volume}{26}
  (\bibinfo{year}{2010}) \bibinfo{pages}{065008}. \URLprefix
  \url{https://doi.org/10.1088/0266-5611/26/6/065008}.
  \DOIprefix\doi{10.1088/0266-5611/26/6/065008}.
\bibitem[{Herman et~al.(2012)Herman, Gardu\~{n}o, Davidi, and Censor}]{hgdc12}
\bibinfo{author}{G.~T. Herman}, \bibinfo{author}{E.~Gardu\~{n}o},
  \bibinfo{author}{R.~Davidi}, \bibinfo{author}{Y.~Censor},
\newblock \bibinfo{title}{Superiorization: An optimization heuristic for
  medical physics},
\newblock \bibinfo{journal}{Medical Physics} \bibinfo{volume}{39}
  (\bibinfo{year}{2012}) \bibinfo{pages}{5532--5546}.
\bibitem[{Censor and Zaslavski(2015)}]{CZ15}
\bibinfo{author}{Y.~Censor}, \bibinfo{author}{A.~J. Zaslavski},
\newblock \bibinfo{title}{Strict {F}ej\'{e}r monotonicity by superiorization of
  feasibility-seeking projection methods},
\newblock \bibinfo{journal}{Journal of Optimization Theory and Applications}
  \bibinfo{volume}{165} (\bibinfo{year}{2015}) \bibinfo{pages}{172--187}.
\bibitem[{Censor et~al.(2021)Censor, Gardu{\~n}o, Helou, and Herman}]{cghh21}
\bibinfo{author}{Y.~Censor}, \bibinfo{author}{E.~Gardu{\~n}o},
  \bibinfo{author}{E.~S. Helou}, \bibinfo{author}{G.~T. Herman},
\newblock \bibinfo{title}{Derivative-free superiorization: {Principle} and
  algorithm},
\newblock \bibinfo{journal}{Numerical Algorithms} \bibinfo{volume}{88}
  (\bibinfo{year}{2021}) \bibinfo{pages}{227--248}.
\bibitem[{Aragón-Artacho et~al.(2023)Aragón-Artacho, Censor, Gibali, and
  Torregrosa-Belén}]{ACGT23}
\bibinfo{author}{F.~J. Aragón-Artacho}, \bibinfo{author}{Y.~Censor},
  \bibinfo{author}{A.~Gibali}, \bibinfo{author}{D.~Torregrosa-Belén},
\newblock \bibinfo{title}{The superiorization method with restarted
  perturbations for split minimization problems with an application to
  radiotherapy treatment planning},
\newblock \bibinfo{journal}{Applied Mathematics and Computation}
  \bibinfo{volume}{440} (\bibinfo{year}{2023}) \bibinfo{pages}{127627}.
  \URLprefix
  \url{https://www.sciencedirect.com/science/article/pii/S0096300322007007}.
  \DOIprefix\doi{https://doi.org/10.1016/j.amc.2022.127627}.
\bibitem[{Humphries et~al.(2020)Humphries, Loreto, Halter, O’Keeffe, and
  Ramirez}]{humphries2020comparison}
\bibinfo{author}{T.~Humphries}, \bibinfo{author}{M.~Loreto},
  \bibinfo{author}{B.~Halter}, \bibinfo{author}{W.~O’Keeffe},
  \bibinfo{author}{L.~Ramirez},
\newblock \bibinfo{title}{Comparison of regularized and superiorized methods
  for tomographic image reconstruction},
\newblock \bibinfo{journal}{J. Appl. Numer. Optim} \bibinfo{volume}{2}
  (\bibinfo{year}{2020}) \bibinfo{pages}{77--99}.
\bibitem[{Nikazad et~al.(2022)Nikazad, Abbasi, Afzalipour, and
  Elfving}]{naae22}
\bibinfo{author}{T.~Nikazad}, \bibinfo{author}{M.~Abbasi},
  \bibinfo{author}{L.~Afzalipour}, \bibinfo{author}{T.~Elfving},
\newblock \bibinfo{title}{A new step size rule for the superiorization method
  and its application in computerized tomography},
\newblock \bibinfo{journal}{Numerical Algorithms} \bibinfo{volume}{90}
  (\bibinfo{year}{2022}) \bibinfo{pages}{1253--1277}.
\bibitem[{Cegielski(2012)}]{C12}
\bibinfo{author}{A.~Cegielski}, \bibinfo{title}{Iterative methods for fixed
  point problems in Hilbert spaces}, volume \bibinfo{volume}{2057},
  \bibinfo{publisher}{Springer}, \bibinfo{year}{2012}.
\bibitem[{Clark et~al.(2013)}]{TCIA}
\bibinfo{author}{K.~Clark}, et~al.,
\newblock \bibinfo{title}{The cancer imaging archive ({TCIA}): maintaining and
  operating a public information repository},
\newblock \bibinfo{journal}{Journal of digital imaging} \bibinfo{volume}{26}
  (\bibinfo{year}{2013}) \bibinfo{pages}{1045--1057}.
\bibitem[{Goldgof et~al.(2015)Goldgof, Hall, Hawkins, Schabath, Stringfield,
  Garcia, Gillies et~al.}]{QIN_LUNG_CT}
\bibinfo{author}{D.~Goldgof}, \bibinfo{author}{L.~Hall},
  \bibinfo{author}{S.~Hawkins}, \bibinfo{author}{M.~Schabath},
  \bibinfo{author}{O.~Stringfield}, \bibinfo{author}{A.~Garcia},
  \bibinfo{author}{R.~Gillies}, et~al., \bibinfo{title}{Data from
  {QIN\_LUNG\_CT}. {The Cancer Imaging Archive}}, \bibinfo{year}{2015}.
  \URLprefix \url{http://doi.org/10.7937/K9/TCIA.2015.NPGZYZBZ}.
\bibitem[{Palenstijn et~al.(2011)Palenstijn, Batenburg, and Sijbers}]{PBS11}
\bibinfo{author}{W.~J. Palenstijn}, \bibinfo{author}{K.~J. Batenburg},
  \bibinfo{author}{J.~Sijbers},
\newblock \bibinfo{title}{Performance improvements for iterative electron
  tomography reconstruction using graphics processing units ({GPUs})},
\newblock \bibinfo{journal}{Journal of structural Biology}
  \bibinfo{volume}{176} (\bibinfo{year}{2011}) \bibinfo{pages}{250--253}.
\bibitem[{Van~Aarle et~al.(2016)Van~Aarle, Palenstijn, Cant, Janssens,
  Bleichrodt, Dabravolski, De~Beenhouwer, Batenburg, and Sijbers}]{VPCJB16}
\bibinfo{author}{W.~Van~Aarle}, \bibinfo{author}{W.~J. Palenstijn},
  \bibinfo{author}{J.~Cant}, \bibinfo{author}{E.~Janssens},
  \bibinfo{author}{F.~Bleichrodt}, \bibinfo{author}{A.~Dabravolski},
  \bibinfo{author}{J.~De~Beenhouwer}, \bibinfo{author}{K.~J. Batenburg},
  \bibinfo{author}{J.~Sijbers},
\newblock \bibinfo{title}{Fast and flexible x-ray tomography using the {ASTRA}
  toolbox},
\newblock \bibinfo{journal}{Optics Express} \bibinfo{volume}{24}
  (\bibinfo{year}{2016}) \bibinfo{pages}{25129--25147}.
\bibitem[{Beister et~al.(2012)Beister, Kolditz, and Kalender}]{BKK12}
\bibinfo{author}{M.~Beister}, \bibinfo{author}{D.~Kolditz},
  \bibinfo{author}{W.~A. Kalender},
\newblock \bibinfo{title}{Iterative reconstruction methods in {X-ray CT}},
\newblock \bibinfo{journal}{Physica Medica} \bibinfo{volume}{28}
  (\bibinfo{year}{2012}) \bibinfo{pages}{94--108}.
\bibitem[{Andersen and Kak(1984)}]{andersen1984simultaneous}
\bibinfo{author}{A.~H. Andersen}, \bibinfo{author}{A.~C. Kak},
\newblock \bibinfo{title}{Simultaneous algebraic reconstruction technique
  ({SART}): a superior implementation of the {ART} algorithm},
\newblock \bibinfo{journal}{Ultrasonic Imaging} \bibinfo{volume}{6}
  (\bibinfo{year}{1984}) \bibinfo{pages}{81--94}.
\bibitem[{Censor and Elfving(2002)}]{CE02}
\bibinfo{author}{Y.~Censor}, \bibinfo{author}{T.~Elfving},
\newblock \bibinfo{title}{Block-iterative algorithms with diagonally scaled
  oblique projections for the linear feasibility problem},
\newblock \bibinfo{journal}{SIAM Journal on Matrix Analysis and Applications}
  \bibinfo{volume}{24} (\bibinfo{year}{2002}) \bibinfo{pages}{40--58}.
\bibitem[{Humphries et~al.(2017)Humphries, Winn, and Faridani}]{hwf17}
\bibinfo{author}{T.~Humphries}, \bibinfo{author}{J.~Winn},
  \bibinfo{author}{A.~Faridani},
\newblock \bibinfo{title}{Superiorized algorithm for reconstruction of {CT}
  images from sparse-view and limited-angle polyenergetic data},
\newblock \bibinfo{journal}{Physics in Medicine \& Biology}
  \bibinfo{volume}{62} (\bibinfo{year}{2017}) \bibinfo{pages}{6762}. \URLprefix
  \url{https://dx.doi.org/10.1088/1361-6560/aa7c2d}.
  \DOIprefix\doi{10.1088/1361-6560/aa7c2d}.
\bibitem[{Wang et~al.(2004)Wang, Bovik, Sheikh, and Simoncelli}]{SSIM}
\bibinfo{author}{Z.~Wang}, \bibinfo{author}{A.~C. Bovik},
  \bibinfo{author}{H.~R. Sheikh}, \bibinfo{author}{E.~P. Simoncelli},
\newblock \bibinfo{title}{Image quality assessment: from error visibility to
  structural similarity},
\newblock \bibinfo{journal}{IEEE Trans. Imag. Proc.} \bibinfo{volume}{13}
  (\bibinfo{year}{2004}) \bibinfo{pages}{600--612}.
\bibitem[{Zhang et~al.(2018)Zhang, Liang, Dong, Xie, and Cao}]{zldxc18}
\bibinfo{author}{Z.~Zhang}, \bibinfo{author}{X.~Liang},
  \bibinfo{author}{X.~Dong}, \bibinfo{author}{Y.~Xie},
  \bibinfo{author}{G.~Cao},
\newblock \bibinfo{title}{A sparse-view {CT} reconstruction method based on
  combination of {DenseNet} and deconvolution},
\newblock \bibinfo{journal}{IEEE Trans. Med. Imag.} \bibinfo{volume}{37}
  (\bibinfo{year}{2018}) \bibinfo{pages}{1407--1417}.
\bibitem[{Han and Ye(2018)}]{hy18}
\bibinfo{author}{Y.~Han}, \bibinfo{author}{J.~C. Ye},
\newblock \bibinfo{title}{Framing {U-Net} via deep convolutional framelets:
  Application to sparse-view {CT}},
\newblock \bibinfo{journal}{IEEE Trans. Med. Imag.} \bibinfo{volume}{37}
  (\bibinfo{year}{2018}) \bibinfo{pages}{1418--1429}.
\bibitem[{Humphries et~al.(2019)Humphries, Si, Coulter, Simms, and
  Xing}]{hscsx19}
\bibinfo{author}{T.~Humphries}, \bibinfo{author}{D.~Si},
  \bibinfo{author}{S.~Coulter}, \bibinfo{author}{M.~Simms},
  \bibinfo{author}{R.~Xing},
\newblock \bibinfo{title}{Comparison of deep learning approaches to low dose
  {CT} using low intensity and sparse view data},
\newblock in: \bibinfo{booktitle}{Medical Imaging 2019: Physics of Medical
  Imaging}, volume \bibinfo{volume}{10948},
  \bibinfo{organization}{International Society for Optics and Photonics},
  \bibinfo{year}{2019}, p. \bibinfo{pages}{109484A}.
\bibitem[{Chen et~al.(2018)Chen, Zhang, Chen, Zhang, Zhang, Sun, Lv, Liao,
  Zhou, and Wang}]{czczz18}
\bibinfo{author}{H.~Chen}, \bibinfo{author}{Y.~Zhang},
  \bibinfo{author}{Y.~Chen}, \bibinfo{author}{J.~Zhang},
  \bibinfo{author}{W.~Zhang}, \bibinfo{author}{H.~Sun},
  \bibinfo{author}{Y.~Lv}, \bibinfo{author}{P.~Liao},
  \bibinfo{author}{J.~Zhou}, \bibinfo{author}{G.~Wang},
\newblock \bibinfo{title}{Learn: Learned experts’ assessment-based
  reconstruction network for sparse-data {CT}},
\newblock \bibinfo{journal}{IEEE Trans. Med. Imag.} \bibinfo{volume}{37}
  (\bibinfo{year}{2018}) \bibinfo{pages}{1333--1347}.
\bibitem[{Adler and {\"O}ktem(2017)}]{AO17}
\bibinfo{author}{J.~Adler}, \bibinfo{author}{O.~{\"O}ktem},
\newblock \bibinfo{title}{Solving ill-posed inverse problems using iterative
  deep neural networks},
\newblock \bibinfo{journal}{Inverse Problems} \bibinfo{volume}{33}
  (\bibinfo{year}{2017}) \bibinfo{pages}{124007}.
\bibitem[{Jia et~al.(2022)Jia, McMichael, Mokarzel, Thompson, Si, and
  Humphries}]{JMMTSH22}
\bibinfo{author}{Y.~Jia}, \bibinfo{author}{N.~McMichael},
  \bibinfo{author}{P.~Mokarzel}, \bibinfo{author}{B.~Thompson},
  \bibinfo{author}{D.~Si}, \bibinfo{author}{T.~Humphries},
\newblock \bibinfo{title}{Superiorization-inspired unrolled {SART} algorithm
  with u-net generated perturbations for sparse-view and limited-angle {CT}
  reconstruction},
\newblock \bibinfo{journal}{Physics in Medicine \& Biology}
  \bibinfo{volume}{67} (\bibinfo{year}{2022}) \bibinfo{pages}{245004}.
  \URLprefix \url{https://dx.doi.org/10.1088/1361-6560/aca513}.
  \DOIprefix\doi{10.1088/1361-6560/aca513}.
\bibitem[{Chang et~al.(2017)Chang, Li, Poczos, Kumar, and
  Sankaranarayanan}]{clpk17}
\bibinfo{author}{J.-H.~R. Chang}, \bibinfo{author}{C.-L. Li},
  \bibinfo{author}{B.~Poczos}, \bibinfo{author}{B.~V. Kumar},
  \bibinfo{author}{A.~C. Sankaranarayanan},
\newblock \bibinfo{title}{One network to solve them all -- solving linear
  inverse problems using deep projection models.},
\newblock in: \bibinfo{booktitle}{ICCV}, \bibinfo{year}{2017}, pp.
  \bibinfo{pages}{5889--5898}.
\bibitem[{Gupta et~al.(2018)Gupta, Jin, Nguyen, McCann, and
  Unser}]{gupta2018cnn}
\bibinfo{author}{H.~Gupta}, \bibinfo{author}{K.~H. Jin}, \bibinfo{author}{H.~Q.
  Nguyen}, \bibinfo{author}{M.~T. McCann}, \bibinfo{author}{M.~Unser},
\newblock \bibinfo{title}{{CNN}-based projected gradient descent for consistent
  {CT} image reconstruction},
\newblock \bibinfo{journal}{IEEE Trans. Med. Imag.} \bibinfo{volume}{37}
  (\bibinfo{year}{2018}) \bibinfo{pages}{1440--1453}.
\bibitem[{Li et~al.(2019)Li, Li, Zhang, Montoya, and Chen}]{LLZMC19}
\bibinfo{author}{Y.~Li}, \bibinfo{author}{K.~Li}, \bibinfo{author}{C.~Zhang},
  \bibinfo{author}{J.~Montoya}, \bibinfo{author}{G.-H. Chen},
\newblock \bibinfo{title}{Learning to reconstruct computed tomography images
  directly from sinogram data under a variety of data acquisition conditions},
\newblock \bibinfo{journal}{IEEE Trans. Med. Imag.} \bibinfo{volume}{38}
  (\bibinfo{year}{2019}) \bibinfo{pages}{2469--2481}.
\bibitem[{Zhang et~al.(2017)Zhang, Zuo, Chen, Meng, and Zhang}]{zzcmz17}
\bibinfo{author}{K.~Zhang}, \bibinfo{author}{W.~Zuo},
  \bibinfo{author}{Y.~Chen}, \bibinfo{author}{D.~Meng},
  \bibinfo{author}{L.~Zhang},
\newblock \bibinfo{title}{Beyond a {Gaussian} denoiser: Residual learning of
  deep {CNN} for image denoising},
\newblock \bibinfo{journal}{IEEE Trans. Imag. Proc.} \bibinfo{volume}{26}
  (\bibinfo{year}{2017}) \bibinfo{pages}{3142--3155}.
\bibitem[{Dabov et~al.(2007)Dabov, Foi, Katkovnik, and
  Egiazarian}]{dabov2007image}
\bibinfo{author}{K.~Dabov}, \bibinfo{author}{A.~Foi},
  \bibinfo{author}{V.~Katkovnik}, \bibinfo{author}{K.~Egiazarian},
\newblock \bibinfo{title}{Image denoising by sparse {3-D} transform-domain
  collaborative filtering},
\newblock \bibinfo{journal}{IEEE Trans. Imag. Proc.} \bibinfo{volume}{16}
  (\bibinfo{year}{2007}) \bibinfo{pages}{2080--2095}.
\bibitem[{M{\"a}kinen et~al.(2020)M{\"a}kinen, Azzari, and Foi}]{MAF20}
\bibinfo{author}{Y.~M{\"a}kinen}, \bibinfo{author}{L.~Azzari},
  \bibinfo{author}{A.~Foi},
\newblock \bibinfo{title}{Collaborative filtering of correlated noise: Exact
  transform-domain variance for improved shrinkage and patch matching},
\newblock \bibinfo{journal}{IEEE Trans. Imag. Proc.} \bibinfo{volume}{29}
  (\bibinfo{year}{2020}) \bibinfo{pages}{8339--8354}.
\bibitem[{Elfving et~al.(2012)Elfving, Hansen, and Nikazad}]{EHN12}
\bibinfo{author}{T.~Elfving}, \bibinfo{author}{P.~C. Hansen},
  \bibinfo{author}{T.~Nikazad},
\newblock \bibinfo{title}{Semiconvergence and relaxation parameters for
  projected {SIRT} algorithms},
\newblock \bibinfo{journal}{SIAM Journal on Scientific Computing}
  \bibinfo{volume}{34} (\bibinfo{year}{2012}) \bibinfo{pages}{A2000--A2017}.

\end{thebibliography}
\end{document}